\documentclass[11pt,a4paper]{article}
\usepackage[utf8]{inputenc}
\usepackage[T1]{fontenc}
\usepackage{amsfonts}
\usepackage{amssymb}
\usepackage{amsthm}
\usepackage{amsmath}
\usepackage{graphicx}
\usepackage{mathrsfs}
\usepackage[left=2cm,right=2cm,top=2cm,bottom=2cm]{geometry}

\newtheorem{theorem}{Theorem}
\newtheorem{example}[theorem]{Example}

\newtheorem{lemma}[theorem]{Lemma}

\newtheorem{remark}[theorem]{Remark}

\begin{document}

\title{Optimal Control of Markov Regime-Switching Stochastic Recursive Utilities}
\author{Liangquan Zhang$^{1}$\thanks{L. Zhang acknowledges the financial support partly by the National Nature
Science Foundation of China (Grant No. 11701040, 11871010 \& 61871058) and
Innovation Foundation of BUPT for Youth (No. 500417024 \& 505018304).
E-mail: xiaoquan51011@163.com.} \quad and \quad
Xun Li$^{2}$\thanks{X. Li acknowledges the financial support partly by Hong Kong
RGC under grants 15209614, 15255416 and 15213218. Email: malixun@polyu.edu.hk} \\
%EndAName
{\small 1. School of Science} \\
{\small Beijing University of Posts and Telecommunications} \\
{\small Beijing 100876, China} \\
%EndAName
{\small 2. Department of Applied Mathematics} \\
{\small The Hong Kong Polytechnic University} \\
{\small Hong Kong, China}}

\maketitle

\begin{abstract}
In this paper, we establish a general stochastic maximum principle for optimal control for systems described by a continuous-time Markov regime-switching stochastic recursive utilities model. The control domain is postulated not to be convex, and the diffusion terms depend on control variables. To this end, we first study a kind of classical forward stochastic optimal control problems. Afterwards, based on previous results, we introduce two groups of new first and second-order adjoint equations. The corresponding variational equations for forward-backward stochastic differential equations are derived. In particular, the generator in the maximum principle contains solutions of second-order adjoint equation which is novel. Some interesting examples are concluded as well.
\end{abstract}



\section{Introduction}

\label{sect:1}

A large number of literatures focus on applications of Markov
regime-switching models in finance as well as stochastic optimal control
since the fundamental feature brings from the Markov regime-switching models
in contrast to the traditional models described via the diffusion processes
from the empirical point of view. The basic idea of regime-switching is to
modulate the model with a continuous-time and finite-state Markov chain where
each state represents a regime of the system or level of the economic
indicator. The regime-switching models have been variously applied in stock
trading, option pricing, portfolio selection, risk measurement, etc
(see \cite{LZA, QZ, ZXCS, Zhu, SKM} references therein).

It is well known that stochastic optimal control is one of the central
themes of modern control sciences. Necessary conditions for the optimal
control of the forward continuous-time stochastic system, which are the
so-called stochastic maximum principle of Pontraygin's type, have been
extensively studied since early 1960s. As a matter of fact, the basic idea
of the stochastic maximum principle is to derive a set of necessary
conditions that must be fulfilled by any optimal control. The stochastic
maximum principle describes that any optimal control must satisfy a system of
forward-backward stochastic differential equations (FBSDEs, for short), called the
optimality system, and minimize a functional, called the \textit{Hamiltonian}. It leads to explicit expressions for the optimal controls in
certain cases. The stochastic maximum principle can be applied to situations
where state processes involve random coefficients and state constraints. Recently, there are some works applying the stochastic maximum principle in finance (see Cadenillas and Karatzas \cite{CK}, see references therein).

When Brownian motion is the unique noise source, Peng \cite{P1} initially studied a
general stochastic maximum principle using the second-order adjoint
equations, which allow the control to enter into both drift and
diffusion coefficients under a nonconvex control domain. For more details
about the stochastic maximum principle theory, interested readers may refer
to Yong and Zhou [9]. Forward-backward stochastic control systems driven by
FBSDEs are widely applied in mathematical economics and mathematical finance. %which includes the usual forward SDEs as a special case.
They mainly emerge in stochastic recursive utility optimization problems and principal-agent
problems. Pardoux and Peng \cite{PP1}
proved the well-posedness for nonlinear BSDE. Duffie and Epstein \cite{DE} introduced the notion of
recursive utilities in continuous-time situation, which is actually a type of backward stochastic differential equation (BSDE, for short)
where the generator $f$ is independent of $z$. El Karoui et al \cite{KPQ} extended the recursive utility to the case where $f$ contains $z$. The term $z$ can be interpreted as an ambiguity aversion term in the market (see Chen and Epstein 2002, \cite{CZE}). In particular, the celebrated
Black-Scholes formula indeed provides an effective way of representing the
option price (which is the solution to a kind of linear BSDE) through the
solution to the Black-Scholes equation (parabolic partial differential
equation actually). Since then, BSDE has been extensively studied and employed
in the areas of applied probability and optimal stochastic controls,
particularly in financial engineering (see \cite{KPQ}).

Peng \cite{P1} derived necessary conditions for optimal control for the
partially coupling case when the control domain is convex. Later, Xu \cite{Xu} studied the nonconvex control domain case and obtained the
corresponding necessary conditions. In the above work, however, assume that the
diffusion term in the forward control system does not include the control variable.
Ji and Zhou applied the Ekeland variational principle to
establish a maximum principle for a partially coupled forward-backward
stochastic control system, while the forward state variable is constrained in a
convex set at the terminal time. Meng \cite{Meng} and Wu \cite{Wu1998} obtained the necessary conditions for optimal control of fully
coupled forward-backward stochastic control systems when the control domain
is convex. Shi and Wu \cite{SW} studied the nonconvex control domain
case and obtained the corresponding necessary conditions under some $G$-monotonic assumptions. However, the control variable is not involved in the diffusion term of the forward equation.

In order to study the backward linear-quadratic optimal control problem,
Kohlmann and Zhou \cite{KZ}, Lim and Zhou \cite{LZ} developed a kind of new method
to handle the problem. The variable $z$ is regarded as a control process and
the terminal condition $y_{T}=h(x_{T}) $ as a constraint, and
then it is possible to apply the Ekeland variational principle to obtain the
maximum principle. Adopting this idea, Yong \cite{Yong2010} and Wu \cite{Wu2013} independently established the maximum principle for the recursive stochastic optimal control problem. Nonetheless, the maximum principle derived by these method involves two
unknown parameters. Hence, some hard questions naturally emerge as follows: What is the
second-order variational equation for the BSDE? How to derive the
second-order adjoint equation since the quadratic form with respect to the
variation of $z$. All of which seem to be extremely complicated. Hu \cite{Hms} overcame the above difficulties by introducing two new adjoint equations. Then, the second-order variational equation for BSDE and the
maximum principle are obtained. The main difference of his variational
equations compared with those in Peng \cite{Peng1990} consists in the term $\langle p(t),\delta \sigma(t) \rangle I_{E_{\varepsilon }}(t) $ in the variation of $z$. Due to the
term $\langle p(t),\delta \sigma(t)
\rangle I_{E_{\varepsilon }}(t) $ in the variation of $z$,
Hu developed a global maximum principle, which is novel and different from
those in Wu \cite{Wu2013}, Yong \cite{Yong2010} and other previous works, to solve Peng's open problem. Hu, Ji and Xue further considered the maximum principle
for fully coupled FBSDEs (see \cite{HJX1, HJX2}).

Inspired by above works, we are more interested in studying
optimal control for FBSDE of Markov-regime diffusion type. Compared with
the literature mentioned above, the novelty of the formulation and the
contribution in this paper should be stated as follows: This work states the
necessary condition for stochastic utilities, the maximum principle of ours
is significantly different from the existing results obtained in \cite{Don, ZEK, SKM} which
describe the sufficient condition. In particular, Sun et al \cite{SGZ}
considered the forward-backward controlled systems and established a sufficient
condition for optimal control. This paper fills a much-needed gap for the necessary part.
Due to the presence of regime-switching, the second-order solution
enters in the generator of BSDE which is new. Tao and Wu \cite{TW} also
consider optimal control for FBSDEs modulated by continuous-time and
finite-state Markov chains with convex control domain (see \cite{LW,Pa}). Recently, a
similar work \cite{ZSX} has been done by Zhang et al for forward
mean field system with jump, but not containing the generator. In this paper,
we adopt the idea from Peng \cite{Peng1990} within our own framework using
the Riesz Representation Theorem to interpret the dual relationship between
the adjoint and variational equation, which is much clearer and more readable.
Besides, we prove a general estimate for BSDE (see Lemma \ref{baest}). Then, we can easily derive some delicate estimates for backward variational equation which is important to establish the variational inequality.

The rest of this paper is organized as follows: After some preliminaries and
notation in the second section, we are devoted the third section to
studying the classical forward stochastic control problems adopting the
form, similar to \cite{Peng1990}. Based on the previous estimates, we
establish a general maximum principle for FBSDEs optimal control problem.
Finally, in Section \ref{Remarks}, we conclude the novelty of this paper and
schedule possible generalizations in future. A proof of technique lemma is
put in the Appendix \ref{app}.

\section{Preliminaries}

\label{sect2}

Throughout this paper, we denote by $\mathbb{R}^{n}$ the space of $n$-dimensional Euclidean space, by $\mathbb{R}^{n\times d}$ the space the
matrices with order $n\times d$. Let $(\Omega,\mathcal{F},P)$ be a completed filtered probability space. For a fixed finite time $T$
horizon, $\{ W_{t},0\leq t\leq T\}$ is an $\mathbb{R}^{d}$-valued standard Brownian motion defined on the space.

To model the controlled state process, we first introduce a set of
Markov jump martingales associated with the chain $\alpha $. We adopt the
notation from Elliott et al \cite{EAM}. The state space of the
chain with a finite state space is denoted by $S:=\{ e_{1},e_{2},\ldots
,e_{D}\} $, where $D\in N$, $e_{i}\in \mathbb{R}^{D}$, and the $j$-th
component of $e_{i}$ is the Kronecker delta $\delta_{ij}$ for each $1\leq
i,j\leq D$. The state space $S$ is called a canonical state space. We
suppose that the chain is homogeneous and irreducible. To specify
statistical or probabilistic properties of the chain $\alpha $, we define
the generator $Q:=[ \lambda_{ij}]_{i,j=1,2,\ldots,D}$ of the
chain under $P$. This is also called the rate matrix, or the $Q$-matrix. For
each $i,j=1,2,\ldots,D$, $\lambda_{ij}$ is the constant transition
intensity of the chain from state $e_{i}$ to state $e_{j}$ at time $t$. Note
that $\lambda_{ij}\geq 0$ for $i\neq j$ an $\sum_{j=1}^{D}\lambda
_{ij}=0 $, so $\lambda_{ii}\leq 0$. In general, for each $i,j=1,2,\ldots,D$
with $i\neq j$, we suppose that $\lambda_{ij}>0$, which means that $\lambda
_{ii}<0$. From Elliott et al \cite{EAM}, we have the following
semimartingale expression for $\alpha$:
\begin{equation*}
\alpha(t) =\alpha(0) +\int_{0}^{t}Q^{\top }\alpha(s) \mathrm{d}s+M(t),
\end{equation*}
where $\{ M(t) \}_{0\leq t\leq T}$ is an $\mathbb{R}%
^{D}$-valued, $(\mathcal{F},P)$-martingale. For each $i,j=1,2,\ldots,D$,
with $i\neq j$, let $\mathcal{J}_{ij}(t) $ be the number of
jumps from state $e_{i}$ to state $e_{j}$ up to time $t$. Then
\begin{eqnarray*}
\mathcal{J}_{ij}(t) &=&\sum_{0<s\leq t}\langle \alpha(s-),e_{i}\rangle \langle \alpha(s),e_{j}\rangle \\
&=&\sum_{0<s\leq t}\langle \alpha(s-),e_{i}\rangle\langle \alpha(s) -\alpha(s-),e_{j}\rangle \\
&=&\int_{0}^{t}\langle \alpha(s-),e_{i}\rangle\langle Q^{\top }\alpha(s),e_{j}\rangle \mathrm{d}s+\int_{0}^{t}\langle \alpha(s-),e_{i}\rangle\langle \mathrm{d}M(s),e_{j}\rangle \mathrm{d}s \\
&=&\lambda_{ij}\int_{0}^{t}\langle \alpha(s-),e_{i}\rangle \mathrm{d}s+m_{ij}(t),
\end{eqnarray*}
where
\begin{equation*}
m_{ij}=\int_{0}^{t}\langle \alpha(s-),e_{i}\rangle
\langle \mathrm{d}M(s),e_{j}\rangle \mathrm{d}s
\end{equation*}
is an $(\mathcal{F},P)$-martingale. The $m_{ij}$'s are called the basic
martingales associated with the chain $\alpha $. For each fixed $%
j=1,2,\ldots,D$, let $\Phi_{j}(t) $ be the number of jumps
into state $e_{j}$ up to time $t$. Thus,
\begin{eqnarray*}
\Phi_{j}(t) &=&\sum_{i=1,i\neq j}^{D}\mathcal{J}_{ij}(t) \\
&=&\sum_{i=1,i\neq j}^{D}\lambda_{ij}\int_{0}^{t}\langle \alpha(s-),e_{i}\rangle \mathrm{d}s+\sum_{i=1,i\neq j}^{D}m_{ij}(t) .
\end{eqnarray*}
Set $\tilde{\Phi}_{j}(t) :=\sum_{i=1,i\neq j}^{D}m_{ij}(
t)$. For each $j=1,2,\ldots,D$, $\{ \tilde{\Phi}_{j}(
t) \}_{0\leq t\leq T}$ is also an $(\mathcal{F},P)$-martingale.
For convenience, define
\begin{equation*}
\lambda_{j}(t) =\sum_{i=1,i\neq j}^{D}\lambda
_{ij}\int_{0}^{t}\langle \alpha(s-),e_{i}\rangle
\mathrm{d}s.
\end{equation*}
Then, for every $j=1,2,\ldots,D$, $\tilde{\Phi}_{j}(t) =\Phi
_{j}(t) -\lambda_{j}(t) $ is also an $(\mathcal{F}%
,P)$-martingale.

Note that $\top $ appearing in this paper as superscript denotes the
transpose of a matrix. Let $U$ be a compact subset of $\mathbb{R}^{k}$. In
what follows, $K$ represents a generic constant, which can be different from
line to line.

The state process is governed by FBSDEs of the following type:
\begin{eqnarray}
x(t) &=&x+\int_{0}^{t}b(s,x(s),u(s),\alpha(s)) \mathrm{d}s+\int_{0}^{t}\sigma(s,x(s),u(s),\alpha(s))\mathrm{d}W(s)  \notag \\
&&+\int_{0}^{t}\gamma(s,x(s-),u(s-),\alpha(s-)) \mathrm{d}\tilde{\Phi}(s),\label{sde1} \\
y(t) &=&g(x(T),\alpha(T))+\int_{t}^{T}f(s,x(s),y(s),z(s),\kappa(s),u(s),\alpha(s))
\mathrm{d}s  \notag \\
&&-\int_{t}^{T}z(s) dW(s) -\int_{t}^{T}\kappa(s) \mathrm{d}\tilde{\Phi}(s),  \label{bsde1}
\end{eqnarray}
for all $t\in [ 0,T]$, $P$-a.s., with $X(t-)=x$, on some
filtered probability space $(\Omega,\mathcal{F},P) $, where
\begin{eqnarray*}
b &:&[ 0,T] \times \mathbb{R}^{L}\times U\times S\rightarrow
\mathbb{R}^{L}, \\
\sigma &:&[ 0,T] \times \mathbb{R}^{L}\times U\times S\rightarrow
\mathbb{R}^{L\times d}, \\
\gamma &:&[ 0,T] \times \mathbb{R}^{L}\times U\times S\rightarrow
\mathbb{R}^{L\times D}, \\
f &:&[ 0,T] \times \mathbb{R}^{L}\times \mathbb{R\times R}^{d}\times \mathbb{R}^{D}\mathbb{\times }U\times S\rightarrow \mathbb{R}^{L\times d}, \\
g &:&\mathbb{R}^{L}\times S\rightarrow \mathbb{R}
\end{eqnarray*}
are given deterministic functions. Denote $\mathcal{U}_{ad}$ as the
admissible control set, where
\begin{equation*}
\mathcal{U}_{ad} = \left\{u(\cdot): [0,T]\times \Omega \rightarrow U
\left\vert\begin{array}{l} u(\cdot)\text{ is an }\mathcal{F}_{t}\text{-adapted process with
values in }U \\
\text{satisfying }\mathbb{E}\Big[\int_{0}^{T}\vert u(s)\vert ^{\beta}\mathrm{d}s\Big] <0\end{array}\right.\right\}.
\end{equation*}
Our problem is to minimize the following cost functional
over $\mathcal{U}_{ad}:$
\begin{equation}
J(u(\cdot)) =y(0) .  \label{costfbsde}
\end{equation}%
Our stochastic recursive optimal control problem is to find an optimal
control $u(\cdot) \in \mathcal{U}_{ad}$ such that $J(\bar{u}(\cdot)) =\inf\limits_{u(\cdot) \in \mathcal{U}_{ad}}y(0)$.

We now introduce the following spaces of processes:
\begin{align*}
\mathcal{S}^{2}(0,T;\mathbb{R}^{L})\triangleq & \bigg\{ \mathbb{R}^{L}\text{-valued }\mathcal{F}_{t}\text{-adapted process }\phi(\cdot)\text{; }\mathbb{E}\bigg[ \sup\limits_{0\leq t\leq T}\vert \phi(t)\vert^{2}\bigg]<\infty \bigg\}, \\
\mathcal{M}^{2}(0,T;\mathbb{R}^{L\times d})\triangleq & \bigg\{ \mathbb{R}^{L\times d}\text{-valued }\mathcal{F}_{t}\text{-adapted process }\varphi(\cdot)\text{; }\mathbb{E}\bigg[ \int_{0}^{T}\vert \varphi(s)\vert^{2}\mbox{\rm d}s\bigg] <\infty \bigg\}, \\
\mathcal{V}^{2}(0,T;\mathbb{R}^{L\times d})& \triangleq \bigg\{ \mathbb{R}^{L\times d}\text{-valued }\mathcal{F}_{t}\text{-adapted process }\varphi(\cdot)\text{; }\mathbb{E}\bigg[ \int_{0}^{T}\sum_{j=1}^{D}\lambda_j\vert \varphi_j(s) \vert^{2}\mbox{\rm d}s\bigg] <\infty \bigg\} .
\end{align*}
Set $\mathcal{N}^{2}(0,T) \triangleq \mathcal{S}^{2}(0,T;\mathbb{R}^{L})\times \mathcal{S}^{2}(0,T;\mathbb{R}^{L})\times \mathcal{M}^{2}(0,T;\mathbb{R}^{L\times d})\times \mathcal{V}^{2}(0,T;\mathbb{R}^{L\times d})$.
%We make the following assumptions.
Assume that
\begin{description}
\item[(A1)] The coefficients $b$ and $\sigma $ satisfy the Lipschitz
condition for $x$, and $f(t,x,y,z,\kappa,u,e_{i}) $ satisfies
the Lipschitz condition as follows:
\begin{eqnarray*}
&&\vert f(t,x,y,z,\kappa,u,e_{i}) -f(t,x',y',z',\kappa',u,e_{i}) \vert \\
&\leq &K\big(\vert x-x'\vert +\vert y-y'\vert +\vert z-z'\vert +\Vert \kappa -\kappa'\Vert\big) .
\end{eqnarray*}
$x,x^{\prime }\in \mathbb{R}^{L},y,y^{\prime },z,z^{\prime }\in \mathbb{R},\kappa,\kappa^{\prime }\in \mathbb{R}^{D}$, for every $e_{i}\in S$.
Moreover,
\begin{eqnarray*}
&&\vert b(t,x,u,e_{i}) \vert +\vert \sigma(t,x,u,e_{i}) \vert +\vert f(t,x,y,z,\kappa,u,e_{i}) \vert \\
&\leq &K\big(1+\vert x\vert +\vert y\vert+\vert z\vert +\Vert \kappa\Vert\big),
\end{eqnarray*}
for all $e_{i}\in S$, where $\Vert \kappa \Vert
=\sum_{j=1}^{D}\vert \kappa_{j}\vert^{2}\lambda_{j}$.

\item[(A2)] Maps $\theta =(x,y,z,\kappa) \rightarrow b(t,x,u,e_{i})$, $\sigma(t,x,u,e_{i})$, $f(t,\theta,u,e_{i}) $ is twice continuously differential, with the (partial) derivatives up to the second order being uniformly bounded, Lipschitz continuous in $(x,y,z,\kappa) $, and continuous in $u\in U$ and
$e_{i}\in S$. The map $g(x,e_{i}) $ is twice differentiable with the derivatives up to the second order being uniformly bounded and uniformly Lipschitz continuous in $e_{i}\in S$.
\end{description}

Clearly, under Assumptions (A1)-(A2), FBSDEs (\ref{sde1}) and (\ref{bsde1}) admit a unique adapted strong solution $(x,y,z,\kappa) \in \mathcal{N}^{2}(0,T)$ (see \cite{CE}). Due to the influence
of Markov regime switching, we need an extension of It\^{o}'s formula with
Markov regime switching in the sequel.

\begin{lemma}\label{Ito} \sl
Assume that an $L$-dimensional process $x(\cdot) $ driven by
\begin{eqnarray*}
\mathrm{d}x(t) &=&b(t,x(t-),u(t),\alpha(t-)) \mathrm{d}s+\sigma(t,x(t-),u(t),\alpha(t-)) \mathrm{d}W(t) \\
&&+\gamma(t,x(t-),u(t-),\alpha(t-)) \mathrm{d}\tilde{\Phi}(t), \quad t\in [0,T]
\end{eqnarray*}
and the function $\varphi(\cdot,\cdot,e_{i}) \in C^{1,2}([ 0,T] \times \mathbb{R}^{L}) $ for every $e_{i}\in S$.
Then
\begin{eqnarray*}
&&\varphi(T,x(T),\alpha(T)) -\varphi(0,x(0),\alpha(0)) \\
&=&\int_{0}^{T}\frac{\partial \varphi(t,x(t-),\alpha(t-)) }{\partial t}+\sum_{m=1}^{L}\frac{\partial \varphi(t,x(t-),\alpha(t-)) }{\partial x_{m}}b_{m}(t,x(t-),u(t),\alpha(t-)) \\
&&+\frac{1}{2}\sum_{m=1}^{L}\sum_{n=1}^{L}\int_{0}^{T}\frac{\partial^{2}\varphi(t,x(t-),\alpha(t-)) }{\partial x_{m}\partial x_{n}}\sum_{l=1}^{d}\sigma_{ml}\sigma_{ml}(t,x(t-),u(t),\alpha(t-)) \\
&&+\sum_{m=1}^{D}\int_{0}^{T}\Big(\varphi(t,x(t-)+\gamma^{m}(t,x(t-),u(t),\alpha(t-)),e_{m}) -\varphi(t,x(t-),\alpha(t-)) \\
&&-\sum_{n=1}^{L}\frac{\partial \varphi(t,x(t-),\alpha(t-)) }{\partial x_{n}}\gamma_{nm}(t,x(t-),u(t),\alpha(t-))\Big)\lambda_{m}(t) \mathrm{d}t \\
&&+\sum_{m=1}^{L}\frac{\partial \varphi(t,x(t-),\alpha(t-)) }{\partial x_{m}}\sum_{n=1}^{d}\sigma_{mn}(t,x(t-),u(t),\alpha(t-))
\mathrm{d}W^{n}(t) \\
&&+\sum_{m=1}^{D}\int_{0}^{T}\Big(\varphi(t,x(t-)+\gamma^{m}(t,x(t-),u(t),\alpha(t-)),e_{m}) -\varphi(t,x(t-),\alpha(t-)) \Big)\mathrm{d}\tilde{\Phi}(t),
\end{eqnarray*}
where $\gamma^{m}$ denotes the $m$-th columns of the matrix $\gamma$.
\end{lemma}

\begin{lemma}\label{qudr} \sl
Let
\begin{equation*}
\mathrm{d}X_{i}(t)=\gamma_{i}(t,\alpha(t)) \mathrm{d}%
\tilde{\Phi}(t),\text{ }i=1,2
\end{equation*}%
be one-dimensional regime switching equations. Then
\begin{eqnarray*}
X_{1}(t) X_{2}(t) &=&X_{1}(0)
X_{2}(0) +\int_{0}^{t}X_{1}(s-) \mathrm{d}%
X_{2}(s) +\int_{0}^{t}X_{2}(s-) \mathrm{d}%
X_{1}(s) \\
&&+\int_{0}^{t}\sum_{j=1}^{D}\gamma_{1j}(s) \gamma_{2j}(
s) \mathrm{d}\Phi_{j}(s) .
\end{eqnarray*}
\end{lemma}

\textbf{Proof.}
Applying It\^{o}'s formula to $Y(t) =X_{1}(t)
X_{2}(t)$, we have
\begin{eqnarray*}
\mathrm{d}(X_{1}(t) X_{2}(t)) &=&\sum_{j=1}^{D}\gamma_{1j}(t) \gamma_{2j}\lambda_{j}(t) \mathrm{d}t \\
&&+\sum_{j=1}^{D}[ \gamma_{1j}(t) \gamma_{2j}(t) +X_{1}(t-) \gamma_{2j}(t) +X_{2}(t-) \gamma_{1j}(t) ] \mathrm{d}\tilde{\Phi}_{j}(t) \\
&=&\sum_{j=1}^{D}\gamma_{1j}(t) \gamma_{2j}(t)[\lambda_{j}(t) \mathrm{d}t+\mathrm{d}\tilde{\Phi}_{j}(t) ] \\
&&+\sum_{j=1}^{D}[ X_{1}(t-) \gamma_{2}(t)
+X_{2}(t-) \gamma_{1}(t) ] \mathrm{d}\tilde{\Phi}(t) \\
&=&\sum_{j=1}^{D}\gamma_{1j}(t) \gamma_{2j}(t)
\mathrm{d}\Phi_{j}(t) +X_{1}(t-) \mathrm{d}X_{2}(t) +X_{2}(t-) \mathrm{d}X_{1}(t),
\end{eqnarray*}
which implies the desired result.
\hfill $\Box$

\vspace{0.25cm}

\begin{remark}
The process
\begin{eqnarray*}
[X_{1},X_{2}](t)
& \triangleq &\int_{0}^{t}\sum_{j=1}^{D}\gamma_{1j}(s) \gamma_{2j}(s) \mathrm{d}\Phi_{j}(s) \\
&=&\int_{0}^{t}\sum_{j=1}^{D}\gamma_{1j}(s) \gamma_{2j}(s) \lambda_{j}(s) \mathrm{d}s+\int_{0}^{t}\sum_{j=1}^{D}\gamma_{1j}(s) \gamma_{2j}(s) \mathrm{d}\tilde{\Phi}_{j}(s)
\end{eqnarray*}
is called the quadratic covariance of $X_{1}$ and $X_{2}$.
\end{remark}

We now present a technical result which will be useful in the next section. Its
proof can be found in Appendix \ref{app}.

\begin{lemma} \sl
\label{baest} Consider the following linear BSDE with Makov regime switching:
\begin{equation}
\left\{\begin{array}{rcl}
-\mathrm{d}y(t) & = & [A(t) y(t)+B(t) z(t) +C(t) \kappa(t)+F(t,\alpha(t-)) ]\mathrm{d}t \\
&  & -z(s) \mathrm{d}W(s) -\kappa(s)\mathrm{d}\tilde{\Phi}(s), \\
y(T) & = & \xi,
\end{array}\right. \label{lbsde}
\end{equation}
where $\xi \in L^{2}(\mathcal{F}_{T})$ with $\mathbb{E}[\vert \xi \vert^{2}] <\infty$, $A,B:[ 0,T]\times \Omega \rightarrow \mathbb{R}$, $C: [ 0,T] \times \Omega
\rightarrow \mathbb{R}^{D}$ and $F: [ 0,T] \times S\times
\Omega \rightarrow \mathbb{R}$ are $\mathcal{F}_{t}$-adapted, and
\begin{equation}
\left\{\begin{array}{l}
\vert A(t) \vert,\text{ }\vert B(t)
\vert,\text{ }\vert C(t) \vert \leq K, \quad \text{a.e. }t\in [ 0,T],\text{ }P\text{-a.s.}, \\
\displaystyle\int_{0}^{T}[ \mathbb{E}\vert F(s,a(s-))\vert^{2k}]^{\frac{1}{2k}}\mathrm{d}s<\infty,
\end{array}\right.\label{lp1}
\end{equation}
for some $k\geq 1$. Then there exists a positive constant $\tilde{K}$ such that
\begin{eqnarray}
&&\sup_{0\leq t\leq T}\mathbb{E}\big[\vert y(s)\vert^{2k}\big] +\mathbb{E}\bigg[\int_{0}^{T}\Big [\vert z(s) \vert^{2}+\sum_{j=1}^{D}\vert \kappa_{j}(s)\vert^{2}\lambda_{j}(s) \Big ]\mathrm{d}s\bigg] \notag \\
&\leq &\tilde{K}\mathbb{E}\bigg[\vert \xi \vert^{2k}+(\int_{0}^{T}[ \mathbb{E}\vert F(s,\alpha(s-)) \vert^{2k}]^{\frac{1}{2k}}\mathrm{d}s)^{2k}\bigg].  \label{newest}
\end{eqnarray}
\end{lemma}

We outline the proof for the convenience of the reader, which can be found in
the Appendix \ref{app}.

\section{Classical optimal control problems}

\label{section3}

The aim of this section is to derive a kind of variational equation (first
and second-order) and corresponding variational inequality. Observe that the
control variable enters in $\sigma $, $\gamma $ and the control domain $U$
is non-convex, both the convex perturbation and the first-order expansion
approach fail. In contrast to Peng \cite{Peng1990}, one needs to find new
adjoint equations because of presence of Markov chains. We now introduce
the spike variation with respect to optimal control $\bar{u}(\cdot) $, more precisely, let $\varepsilon >0$
and $E_{\varepsilon }\subset [ 0,T] $ be a Borel set with Borel measure $\vert E_{\varepsilon }\vert =\varepsilon$, defined as follows:
\begin{equation*}
u^{\varepsilon }(t) = \left\{\begin{array}{ll}
v, & \text{if }\tau \leq t\leq \tau +\varepsilon, \\
\bar{u}(t), & \text{otherwise,}
\end{array}\right.
\end{equation*}
where $0 \leq \tau <T$ is an arbitrarily fixed time, $\varepsilon >0$ is a
sufficiently small constant, and $v$ is an arbitrary $U$-valued $\mathcal{F}_{\tau }$-measurable random variable such that
$\mathbb{E}\vert v\vert^{3}<+\infty $. Let $x^{\varepsilon }(\cdot) $ be
the trajectory of the control system (\ref{sde1}) corresponding to the control $u^{\varepsilon }(\cdot)$.

The cost functional is
\begin{equation*}
\mathcal{J}(u(\cdot)) =\mathbb{E}\bigg[\int_{0}^{T}l(t,x(t),u(t),\alpha(t-)) \mathrm{d}t+h(x(T),\alpha(T))\bigg],
\end{equation*}
where $l:[ 0,T] \times \mathbb{R}^{L}\times \mathbb{R}^{k}\times I\mathbb{\rightarrow R}$, $h:\mathbb{R}^{L}\mathbb{\rightarrow R}$. Our problem is to find an optimal control
$u(\cdot) \in \mathcal{U}_{ad}$ such that
\begin{equation}
\mathcal{J}(\bar{u}(\cdot)) =\inf_{u(\cdot) \in \mathcal{U}_{ad}}\mathcal{J}(u(\cdot)).  \label{cpr}
\end{equation}
In what follows, we add the following additional assumption:

\begin{description}
\item[(A3)] The map $l$ and $h$ are twice differentiable with the
derivatives up to the second order being uniformly bounded and uniformly
Lipschitz continuous in $u\in U$, $e_{i}\in S$.
\end{description}

We shall derive the variational inequality from the following fact:
$$
J(u^{\varepsilon }(\cdot)) -J(\bar{u}(\cdot)) \geq 0.
$$
Next, fixed $\alpha(\cdot)$, for $\varphi =b,\sigma,\gamma$, we define
\begin{eqnarray*}
\varphi(t) &\triangleq & \varphi(t,\bar{x}(t),\bar{u}(t),\alpha(t)), \\
\varphi_{x}(t) & \triangleq & \varphi_{x}(t,\bar{x}(t),\bar{u}(t),\alpha(t)), \\
\varphi_{xx}(t) &\triangleq &\varphi_{xx}(t,\bar{x}(t),\bar{u}(t),\alpha(t)), \\
\delta \varphi^{\varepsilon }(t) &\triangleq &\varphi(t,\bar{x}(t),u^{\varepsilon }(t),\alpha(t)) -\varphi(t,\bar{x}(t),\bar{u}(t),\alpha(t)), \\
\delta \varphi_{x}^{\varepsilon }(t) &\triangleq &\varphi_{x}(t,\bar{x}(t),u^{\varepsilon }(t),\alpha(t)) -\varphi_{x}(t,\bar{x}(t),\bar{u}(t),\alpha(t)), \\
\delta \varphi_{xx}^{\varepsilon }(t) &\triangleq &\varphi_{xx}(t,\bar{x}(t),u^{\varepsilon }(t),\alpha(t)) -\varphi_{xx}(t,\bar{x}(t),\bar{u}(t),\alpha(t)).
\end{eqnarray*}
Let $x_{1}(\cdot) $ and $x_{2}(\cdot) $ be
respectively solutions to the following stochastic differential equations:
\begin{equation}
\left\{\begin{array}{lll}
\mathrm{d}x_{1}(t) & = & [ b_{x}(t) x_{1}(t) +\delta b^{\varepsilon }(t) ] \mathrm{d}t+\displaystyle\sum_{j=1}^{d}[ \sigma_{x}^{j}(t) x_{1}(t)
+\delta \sigma^{\varepsilon j}(t) ] \mathrm{d}W^{j}(t) \\
&  & +\displaystyle\sum_{j=1}^{D}[ \gamma_{x}^{j}(t) x_{1}(t) +\delta \gamma^{\varepsilon j}(t) ] \mathrm{d}\tilde{\Phi}^{j}(t), \\
x_{1}(0) & = & 0
\end{array}\right.  \label{v1}
\end{equation}
and
\begin{equation}
\left\{\begin{array}{lll}
\mathrm{d}x_{2}(t) & = & [ b_{x}(t) x_{2}(t) +\frac{1}{2}b_{xx}(t) x_{1}(t) x_{1}(t)] \mathrm{d}t \\
&  & +\displaystyle\sum_{j=1}^{d}[ \sigma_{x}^{j}(t) x_{2}(t) +\frac{1}{2}\sigma_{xx}^{j}(t) x_{1}(t)x_{1}(t) ] \mathrm{d}W^{j}(t) \\
&  & +\displaystyle\sum_{j=1}^{D}[ \gamma_{x}^{j}(t) x_{2}(t) +\frac{1}{2}\gamma_{xx}^{j}(t) x_{1}(t)x_{1}(t) ] \mathrm{d}\tilde{\Phi}^{j}(t) \\
&  & +\delta b_{x}^{\varepsilon }(t) x_{1}(t)\mathrm{d}t+\displaystyle\sum_{j=1}^{d}\delta \sigma_{x}^{\varepsilon j}(t)x_{1}(t) \mathrm{d}W^{j}(t) \\
&  & +\displaystyle\sum_{j=1}^{D}\delta \gamma_{x}^{\varepsilon j}(t)x_{1}(t) \mathrm{d}\tilde{\Phi}^{j}(t), \\
x_{2}(0) & = & 0,
\end{array}\right.  \label{v2}
\end{equation}
where
\begin{eqnarray*}
b_{xx}(t) x_{1}(t) x_{1}(t) &\triangleq
& \left(\begin{array}{c}
\text{tr}\{ b_{xx}^{1}(t) \} x_{1}(t)x_{1}^{\top }(t) \\
\vdots \\
\text{tr}\{ b_{xx}^{L}(t) \} x_{1}(t)x_{1}^{\top }(t)
\end{array}\right), \\
\sigma_{xx}^{j}(t) x_{1}(t) x_{1}(t)
& \triangleq & \left(\begin{array}{c}
\text{tr}\{ \sigma_{xx}^{1j}(t) x_{1}(t)
x_{1}^{\top }(t) \} \\
\vdots \\
\text{tr}\{ \sigma_{xx}^{Lj}(t) x_{1}(t)x_{1}^{\top }(t) \}
\end{array}\right),  \quad 1\leq j\leq d
\end{eqnarray*}
and
\begin{equation*}
\gamma_{xx}^{j}(t) x_{1}(t) x_{1}(t)
\triangleq \left(\begin{array}{c}
\text{tr}\{ \gamma_{xx}^{1j}(t) x_{1}(t)x_{1}^{\top }(t) \} \\
\vdots \\
\text{tr}\{ \gamma_{xx}^{Lj}(t) x_{1}(t)x_{1}^{\top }(t)\}
\end{array}\right),   \quad 1\leq j\leq D.
\end{equation*}

\vspace{0.25cm}

\begin{remark}
Equations (\ref{v1}) and (\ref{v2}) are usually called the first and second-order
variational equations, respectively.
\end{remark}

Define the Hamiltonian function as follows:
\begin{eqnarray}
H(t,x,u,e_{i},p,q,s) &\triangleq &l(t,x,u,e_{i})+b^{\top }(t,x,u,e_{i}) p+\text{tr}\big(\sigma(t,x,u,e_{i})q\big)  \notag \\
&&+\sum_{m=1}^{D}\sum_{n=1}^{L}\gamma_{nm}(t,x,u,e_{i})
s_{nm}\lambda_{im},  \label{Hal}
\end{eqnarray}
where $(t,x,u,e_{i},p,q,s) \in [ 0,T] \times \mathbb{R}^{L}\times U\times \mathbb{R}^{L}\mathbb{\times R}^{L\times d}\mathbb{\times R}^{L\times D}$.

The following result proceeds the Taylor expansion of the state with respect
to the control perturbation. We need the following estimations:

\begin{lemma}\label{est} \sl
Let \emph{(A1)-(A2)} hold. Then, we have
\begin{eqnarray}
\sup_{0\leq t\leq T}\mathbb{E}\big[\big\vert x^{\varepsilon }(t) -\bar{x}(t)\big\vert^{2}\big] &=&O(\varepsilon),  \label{e1} \\
\sup_{0\leq t\leq T}\mathbb{E}\big[\big\vert x_{1}(t)\big\vert^{2}\big] &=&O(\varepsilon),  \label{e2} \\
\sup_{0\leq t\leq T}\mathbb{E}\big[\big\vert x_{2}(t)\big\vert^{2}\big] &=&O(\varepsilon^{2}),  \label{e3} \\
\sup_{0\leq t\leq T}\mathbb{E}\big[\big\vert x^{\varepsilon }(t) -\bar{x}(t) -x_{1}(t)\big\vert^{2}\big] &=&O(\varepsilon^{2}),  \label{e4} \\
\sup_{0\leq t\leq T}\mathbb{E}\big[\big\vert x^{\varepsilon }(t) -\bar{x}(t) -x_{1}(t) -x_{2}(t)\big\vert^{2}\big] &=&o(\varepsilon^{2}).  \label{e5}
\end{eqnarray}
\end{lemma}

\paragraph{Proof.}

Let $\xi^{\varepsilon}(t) \triangleq x^{\varepsilon}(t) -\bar{x}(t)$, $t\in [0,T] $. Then, we have
\begin{equation}
\left\{\begin{array}{lll}
\mathrm{d}\xi^{\varepsilon}(t)
& = &(\hat{b}^{\varepsilon}(t) \xi^{\varepsilon}(t) +\delta b(t)) \mathrm{d}t+(\hat{\sigma}^{\varepsilon} t) \xi^{\varepsilon}(t) +\delta\sigma(t)) \mathrm{d}W(t)  \\
&  & +(\hat{\gamma}^{\varepsilon }(t)\xi^{\varepsilon}(t) +\delta\gamma(t)) \mathrm{d}\tilde{\Phi}(t), \\
\xi^{\varepsilon}(t)  & = & 0,
\end{array}\right.   \label{o1}
\end{equation}
where
\begin{eqnarray*}
\hat{b}^{\varepsilon }(t)  & \triangleq & \int_{0}^{1}b_{x}(t,\bar{x}(t+\lambda(x^{\varepsilon }(t) -\bar{x}(t))),u^{\varepsilon }(t),\alpha(t-))\mathrm{d}\lambda, \\
\hat{\sigma}^{\varepsilon }(t)  & \triangleq & \int_{0}^{1}\sigma_{x}(t,\bar{x}(t+\lambda(x^{\varepsilon }(t) -\bar{x}(t))),u^{\varepsilon }(t),\alpha(t-)) \mathrm{d}\lambda, \\
\hat{\gamma}^{\varepsilon }(t)  & \triangleq & \int_{0}^{1}\gamma_{x}(t,\bar{x}(t+\lambda(x^{\varepsilon }(t) -\bar{x}(t))),u^{\varepsilon }(t),\alpha(t-))\mathrm{d}\lambda .
\end{eqnarray*}
By the classical method, it is easy to verify $\sup\limits_{0\leq t\leq T}\mathbb{E}
[|x_{1}(t)|^{2p}]=O\left( \varepsilon ^{p}\right) $ and $\sup\limits_{0\leq t\leq T}\mathbb{E}[|x_{2}(t)|^{2p}]=O\left( \varepsilon ^{2p}\right) $ for $p\geq 1.$ Taking the expectation after squaring both sides of (\ref{o1}) and using Gronwall inequality, we have
\begin{eqnarray*}
&&\sup_{0\leq t\leq T}\mathbb{E}[\vert \xi^{\varepsilon }(t)\vert^{2}]  \\
&\leq &C\mathbb{E}\bigg[\bigg\vert \int_{0}^{T}\hat{b}^{\varepsilon}(s) \xi^{\varepsilon }(s) +\delta b(s) \mathrm{d}s\bigg\vert^{2} \\
&&+\bigg\vert \int_{0}^{T}\hat{\sigma}^{\varepsilon }(s) \xi^{\varepsilon }(s) +\delta \sigma(s) \mathrm{d}W(s)\bigg\vert^{2}+\bigg\vert \int_{0}^{T}\hat{\gamma}^{\varepsilon }(s) \xi^{\varepsilon }(s) +\delta\gamma(s) \mathrm{d}\tilde{\Phi}(s)\bigg\vert^{2}\bigg] \\
&\leq &C \mathbb{E}\bigg [\int_{0}^{T}I_{E_{\varepsilon }}(s)
\mathrm{d}s\int_{0}^{T}\vert \hat{b}^{\varepsilon }(s) \xi^{\varepsilon }(s) +\delta b(s) \vert^{2}\mathrm{d}t+\int_{0}^{T}\vert \hat{\sigma}^{\varepsilon}(s) \xi^{\varepsilon }(s) +\delta\sigma(s)\vert^{2}\mathrm{d}s \\
&&+\int_{0}^{T}\Big\vert\sum_{j=1}^{D}\sum_{m=1}^{L}(\hat{\gamma}^{\varepsilon j}(s) \xi^{\varepsilon }(s) +\delta\gamma^{j}(s))_{m}\Big\vert^{2}\lambda_{j}(s)\mathrm{d}s\bigg] \\
&\leq &C\varepsilon.
\end{eqnarray*}
Similarly, we can prove (\ref{e2}). We now deal with (\ref{e3}). In fact,
\begin{eqnarray*}
&&\sup_{0\leq t\leq T}\mathbb{E}[ \vert x_{2}(t)\vert^{2}]  \\
&\leq &C\mathbb{E}\bigg [\Big |\int_{0}^{T}b_{x}(t) x_{2}(t) +\frac{1}{2}b_{xx}(t) x_{1}(t) x_{1}(t) \mathrm{d}t\Big\vert^{2} \\
&&+\int_{0}^{T}I_{E_{\varepsilon }}(s) \mathrm{d}s\int_{0}^{T}\Big |\sum_{j=1}^{d}[ \sigma_{x}^{j}(t) x_{2}(t) +\frac{1}{2}\sigma_{xx}^{j}(t) x_{1}(t)x_{1}(t) ] \Big\vert^{2}\mathrm{d}t \\
&&+\int_{0}^{T}I_{E_{\varepsilon }}(s) \mathrm{d}s\int_{0}^{T}\Big |\sum_{j=1}^{D}[ \gamma_{x}^{j}(t) x_{2}(t) +\frac{1}{2}\gamma_{xx}^{j}(t) x_{1}(t)x_{1}(t)] \lambda_{j}(t) \Big\vert^{2}\mathrm{d}t \\
&&+\Big |\int_{0}^{T}\delta b_{x}^{\varepsilon }(t) x_{1}(t) \mathrm{d}t\Big\vert^{2}+\int_{0}^{T}I_{E_{\varepsilon }}(s) \mathrm{d}s\int_{0}^{T}\Big |\sum_{j=1}^{d}\delta \sigma_{x}^{\varepsilon j}(t) x_{1}(t) \Big\vert^{2}\mathrm{d}t \\
&&+\int_{0}^{T}I_{E_{\varepsilon }}(s) \mathrm{d}s\int_{0}^{T}\Big |\sum_{j=1}^{D}\delta \gamma_{x}^{\varepsilon j}(t)x_{1}(t) \lambda_{j}(t) \Big\vert^{2}\mathrm{d}t\bigg] \\
&\leq &C\varepsilon^{2}.
\end{eqnarray*}
To prove (\ref{e4}), let
\begin{eqnarray*}
\eta^{\varepsilon }(t)  &=&x^{\varepsilon }(t) -\bar{x}(t) -x_{1}(t)  \\
&=&\xi^{\varepsilon }(t) -x_{1}(t) .
\end{eqnarray*}
Then, it follows from (\ref{o1}) and (\ref{v1}) that we obtain
\begin{eqnarray*}
\mathrm{d}\eta^{\varepsilon }(t)  &=&\big[\hat{b}^{\varepsilon }(t) \eta^{\varepsilon }(t) +(\hat{b}^{\varepsilon }(t) -b_{x}(t))x_{1}(t) +\delta b(t) -\delta b^{\varepsilon }(t)\big]\mathrm{d}t \\
&&+\big[\hat{\sigma}^{\varepsilon }(t) \xi^{\varepsilon }(t) +\delta\sigma(t) -\sigma_{x}(t)x_{1}(t) -\delta\sigma^{\varepsilon}(t)\big]\mathrm{d}W(t)  \\
&&+\big[\hat{\gamma}^{\varepsilon }(t)\xi^{\varepsilon}(t) +\delta \gamma(t) -\gamma_{x}(t)x_{1}(t) -\delta\gamma^{\varepsilon}(t)\big]\mathrm{d}\tilde{\Phi}(t)  \\
&=& \big[\hat{b}^{\varepsilon }(t) \eta^{\varepsilon }(t) +(\hat{b}^{\varepsilon }(t) -b_{x}(t)) x_{1}(t) +\delta b(t) -\delta b^{\varepsilon}(t)\big]\mathrm{d}t \\
&&+\big[\hat{\sigma}^{\varepsilon }(t)\eta^{\varepsilon}(t) +(\hat{\sigma}^{\varepsilon }(t) -\sigma_{x}(t)) x_{1}(t) +\delta \sigma(t) -\delta \sigma^{\varepsilon }(t)\big]\mathrm{d}W(t)  \\
&&+\big[(\hat{\gamma}^{\varepsilon }(t) \eta^{\varepsilon }(t) +(\hat{\gamma}^{\varepsilon }(t) -\gamma_{x}(t))x_{1}(t) +\delta\gamma(t) -\delta \gamma^{\varepsilon }(t))\big]\mathrm{d}\tilde{\Phi}(t) .
\end{eqnarray*}
Therefore,
\begin{eqnarray*}
& & \sup_{0\leq t\leq T}\mathbb{E}\big[\vert \eta^{\varepsilon}(t)\vert^{2}\big] \\
&\leq &\mathbb{E}\bigg\vert \int_{0}^{T}[ \hat{b}^{\varepsilon}(t) \eta^{\varepsilon }(t) +(\hat{b}^{\varepsilon}(t) -b_{x}(t)) x_{1}(t) +\delta
b(t) -\delta b^{\varepsilon }(t) ]\mathrm{d}t\bigg\vert^{2} \\
&&+\mathbb{E}\int_{0}^{T}\big\vert[\hat{\sigma}^{\varepsilon}(t)\eta^{\varepsilon }(t) +(\hat{\sigma}^{\varepsilon }(t) -\sigma_{x}(t))
x_{1}(t) +\delta\sigma(t) -\delta\sigma^{\varepsilon }(t)]\big\vert^{2}\mathrm{d}t \\
&&+\int_{0}^{T}\sum_{j=1}^{D}\sum_{m=1}^{L}\bigg\vert\Big[\hat{\gamma}^{\varepsilon j}(t)\eta^{\varepsilon }(t) +(
\hat{\gamma}^{\varepsilon j}(t) -\gamma_{x}^{j}(t)) x_{1}(t)  +\delta \gamma^{j}(t) -\delta \gamma^{\varepsilon j}(
t) \Big]_{m}\bigg\vert^{2}\lambda_{j}(s) \mathrm{d}s.
\end{eqnarray*}
We consider the following term
\begin{eqnarray*}
& & \int_{0}^{T}\vert \hat{b}^{\varepsilon }(t) -b_{x}(t) \vert^{2}\vert x_{1}(t) \vert^{2}\mathrm{d}t \\
&\leq &\varepsilon \int_{0}^{T}\vert \int_{0}^{1}[b_{x}(t,\bar{x}(t+\lambda(\vert x^{\varepsilon }(t) -\bar{x}(t)\vert)),u^{\varepsilon}(t),\alpha(t-)) -b_{x}(t)]\mathrm{d}\lambda \vert^{2}\mathrm{d}t \\
&\leq &\varepsilon C\int_{0}^{T}(\vert x^{\varepsilon }(t) -\bar{x}(t) \vert^{2}+\vert\delta b_{x}(t) \vert^{2}) \mathrm{d}t \\
&\leq &C\varepsilon^{2}.
\end{eqnarray*}
Similarly, we have
\begin{equation*}
\int_{0}^{T}\vert \hat{\sigma}^{\varepsilon }(t) -\sigma_{x}(t) \vert^{2}\vert x_{1}(t)\vert^{2}\mathrm{d}t\leq C\varepsilon^{2}
\end{equation*}
and
\begin{equation*}
\int_{0}^{T}\vert \hat{\gamma}^{\varepsilon j}(t) -\gamma_{x}^{j}(t) \vert^{2}\vert x_{1}(t)\vert^{2}\mathrm{d}t\leq C\varepsilon^{2}.
\end{equation*}
Hence, we get (\ref{e4}). Set
\begin{equation*}
x_{3}(t) =x_{1}(t) +x_{2}(t) .
\end{equation*}
Then, we have
\begin{eqnarray*}
& & \int_{0}^{t}b(s,\bar{x}(s) +x_{3}(s),u^{\varepsilon }(s),\alpha(s-)) \mathrm{d}s \\
& & +\int_{0}^{t}\sigma(s,\bar{x}(s) +x_{3}(s),u^{\varepsilon }(s),\alpha(s-)) \mathrm{d}W(s)  \\
& & +\int_{0}^{t}\gamma(s,\bar{x}(s) +x_{3}(s),u^{\varepsilon }(s),\alpha(s-)) \mathrm{d}\tilde{\Phi}(s)  \\
&=&\int_{0}^{t}b(s,\bar{x}(s),u^{\varepsilon }(s),\alpha(s-)) +b_{x}(s,\bar{x}(s),u^{\varepsilon }(s),\alpha(s-))x_{3}(s) \mathrm{d}s \\
&&+\int_{0}^{1}\int_{0}^{1}\lambda b_{xx}(s,\bar{x}(s)+\lambda \mu x_{3}(s),u^{\varepsilon }(s),\alpha(s-)) \mathrm{d}\lambda \mathrm{d}\mu x_{3}(s) x_{3}(s) \mathrm{d}s \\
&&+\int_{0}^{t}\sigma(s,\bar{x}(s),u^{\varepsilon}(s),\alpha(s-)) +\sigma_{x}(s,\bar{x}(s),u^{\varepsilon}(s),\alpha(s-)) x_{3}(s)\mathrm{d}W(s)  \\
&&+\int_{0}^{1}\int_{0}^{1}\lambda \sigma_{xx}(s,\bar{x}(s) +\lambda \mu x_{3}(s),u^{\varepsilon}(s),\alpha(s-))\mathrm{d}\lambda \mathrm{d}\mu x_{3}(s) x_{3}(s) \mathrm{d}W(s)  \\
&&+\int_{0}^{t}\gamma(s,\bar{x}(s),u^{\varepsilon}(s),\alpha(s-)) +\gamma_{x}(s,\bar{x}(s),u^{\varepsilon }(s),\alpha(s-)) x_{3}(s) \mathrm{d}\tilde{\Phi}(s)  \\
&&+\int_{0}^{1}\int_{0}^{1}\lambda\gamma_{xx}(s,\bar{x}(s) +\lambda \mu x_{3}(s),u^{\varepsilon }(s),\alpha(s-))\mathrm{d}\lambda \mathrm{d}\mu x_{3}(s)x_{3}(s)\mathrm{d}\tilde{\Phi}(s)  \\
&=&\int_{0}^{t}b(s) \mathrm{d}s+\int_{0}^{t}\sigma(s)\mathrm{d}W(s) +\int_{0}^{t}b_{x}(s)\mathrm{d}s+\int_{0}^{t}\sigma_{x}(s)\mathrm{d}W(s)  \\
&&+\int_{0}^{t}(b(s,\bar{x}(s),u^{\varepsilon}(s),\alpha(s-)) -b(s))\mathrm{d}s \\
&&+\int_{0}^{t}(\sigma(s,\bar{x}(s),u^{\varepsilon}(s),\alpha(s-)) -\sigma(s))\mathrm{d}W(s)  \\
&&+\frac{1}{2}\int_{0}^{t}b_{xx}(s) x_{3}(s)x_{3}(s)\mathrm{d}s+\frac{1}{2}\int_{0}^{t}\sigma_{xx}(s) x_{3}(s) x_{3}(s)\mathrm{d}W(s)
\end{eqnarray*}
\begin{eqnarray}
&&+\int_{0}^{t}(b_{x}(s,\bar{x}(s),u^{\varepsilon}(s),\alpha(s-)) -b_{x}(s))\mathrm{d}s \notag \\
&&+\int_{0}^{t}(\sigma_{x}(s,\bar{x}(s),u^{\varepsilon }(s),\alpha(s-)) -\sigma_{x}(s)) \mathrm{d}W(s)   \notag \\
&&+\int_{0}^{t}\int_{0}^{1}\int_{0}^{1}\lambda \big [b_{xx}(s,\bar{x}(s) +\lambda \mu x_{3}(s),u^{\varepsilon }(s),\alpha(s-))-b_{xx}(s) \big ]\mathrm{d}\lambda \mathrm{d}\mu~ x_{3}(s) x_{3}(s) \mathrm{d}s  \notag \\
&&+\int_{0}^{t}\int_{0}^{1}\int_{0}^{1}\lambda \big[\sigma_{xx}(s,\bar{x}(s) +\lambda \mu x_{3}(s),u^{\varepsilon}(s),\alpha(s-))-\sigma_{xx}(s) \big ]\mathrm{d}\lambda \mathrm{d}\mu ~ x_{3}(s) x_{3}(s) \mathrm{d}W(s)   \notag \\
&&+\int_{0}^{t}\int_{0}^{1}\int_{0}^{1}\lambda \big [\gamma_{xx}(s,\bar{x}(s) +\lambda \mu x_{3}(s),u^{\varepsilon}(s),\alpha(s-))-\gamma_{xx}(s) \big ]\mathrm{d}\lambda \mathrm{d}\mu ~ x_{3}(s) x_{3}(s) \mathrm{d}\tilde{\Phi}(s)   \notag \\
&=&\bar{x}(t) +x_{3}(t) -x+\int_{0}^{t}\Pi^{\varepsilon }(s)\mathrm{d}s+\int_{0}^{t}\Theta^{\varepsilon}(s) \mathrm{d}W(s) +\int_{0}^{t}\Xi^{\varepsilon}(s) \mathrm{d}\tilde{\Phi}(s),  \label{l1}
\end{eqnarray}
\begin{eqnarray*}
\Pi^{\varepsilon }(s)  &=&\frac{1}{2}b_{xx}(s)(x_{2}(s) x_{2}(s) +2x_{1}(s)x_{2}(s))  \\
&&+(b_{x}(s,\bar{x}(s),u^{\varepsilon }(s),\alpha(s-)) -b_{x}(s))x_{2}(s)  \\
&&+\int_{0}^{1}\int_{0}^{1}\lambda b_{xx}(s,\bar{x}(s)+\lambda \mu x_{3}(s),u^{\varepsilon }(s),\alpha(s-)) -b_{xx}(s) \mathrm{d}\lambda \mathrm{d}\mu ~ x_{3}(s)x_{3}(s),
\end{eqnarray*}
\begin{eqnarray*}
\Theta^{\varepsilon }(s)  &=&\frac{1}{2}\sigma_{xx}(s)(x_{2}(s) x_{2}(s) +2x_{1}(s) x_{2}(s))  \\
&&+(\sigma_{x}(s,\bar{x}(s),u^{\varepsilon }(s),\alpha(s-)) -\sigma_{x}(s)) x_{2}(s)  \\
&&+\int_{0}^{1}\int_{0}^{1}\lambda \sigma_{xx}(s,\bar{x}(s) +\lambda \mu x_{3}(s),u^{\varepsilon }(s),\alpha(s-)) -\sigma_{xx}(s) \mathrm{d}\lambda \mathrm{d}\mu  ~ x_{3}(s) x_{3}(s)
\end{eqnarray*}
and
\begin{eqnarray*}
\Xi^{\varepsilon }(s)  &=&\frac{1}{2}\gamma_{xx}(s)(x_{2}(s) x_{2}(s) +2x_{1}(s) x_{2}(s))  \\
&&+(\gamma_{x}(s,\bar{x}(s),u^{\varepsilon }(s),\alpha(s-)) -\gamma_{x}(s)) x_{2}(s)  \\
&&+\int_{0}^{1}\int_{0}^{1}\lambda \gamma_{xx}(s,\bar{x}(s) +\lambda \mu x_{3}(s),u^{\varepsilon }(s),\alpha(s-)) -\gamma_{xx}(s) \mathrm{d}\lambda \mathrm{d}\mu  ~ x_{3}(s) x_{3}(s).
\end{eqnarray*}
One can check that
\begin{equation}
\sup_{0\leq t\leq T}\mathbb{E}\bigg[\bigg\vert \int_{0}^{t}\Pi^{\varepsilon }(s) \mathrm{d}s\bigg\vert^{2}+\bigg\vert\int_{0}^{t}\Theta^{\varepsilon }(s) \mathrm{d}W(s)\bigg\vert^{2}+\bigg\vert \int_{0}^{t}\Xi^{\varepsilon }(s)\mathrm{d}\tilde{\Phi}(s)\bigg\vert^{2}\bigg] =o(\varepsilon^{2}).  \label{e6}
\end{equation}
Recall
\begin{eqnarray*}
x^{\varepsilon }(t)  &=&x+\int_{0}^{t}b(s,x^{\varepsilon}(s),u^{\varepsilon }(s),\alpha(s-)) \mathrm{d}s+\int_{0}^{t}\sigma(s,x^{\varepsilon }(s),u^{\varepsilon }(s),\alpha(s-))\mathrm{d}W(s)  \\
&&+\int_{0}^{t}\gamma(s,x^{\varepsilon }(s),u^{\varepsilon }(s),\alpha(s-)) \mathrm{d}\tilde{\Phi}(s).
\end{eqnarray*}
From (\ref{l1}), we have
\begin{eqnarray*}
x^{\varepsilon }(t) -\bar{x}(t) -x_{3}(t)  &=&\int_{0}^{t}A^{\varepsilon }(s)(x^{\varepsilon }(s) -\bar{x}(s) -x_{3}(s))\mathrm{d}s \\
&&+\int_{0}^{t}B^{\varepsilon }(s)(x^{\varepsilon }(s) -\bar{x}(s) -x_{3}(s)) \mathrm{d}W(s)  \\
&&+\int_{0}^{t}\Pi^{\varepsilon }(s)\mathrm{d}s+\int_{0}^{t}\Theta^{\varepsilon }(s) \mathrm{d}W(s) +\int_{0}^{t}\Xi^{\varepsilon }(s)\mathrm{d}\tilde{\Phi}(s) .
\end{eqnarray*}
It is easy to check that
\begin{equation*}
\vert A^{\varepsilon }(s,\omega)\vert +\vert B^{\varepsilon }(s,\omega)\vert \leq C,  \quad 0\leq s\leq T,  \quad \forall \omega \in \Omega .
\end{equation*}
By virtue of Gronwall inequality, It\^{o}'s formula and (\ref{e6}), we complete the proof.
\hfill $\Box $

\begin{lemma}\label{varineq} \sl
Let \emph{(A1)-(A3)} hold. Then
\begin{eqnarray}
&&\mathbb{E}\bigg[\int_{0}^{T}l_{x}(s,\bar{x}(s),\bar{u}(s))(x_{1}(s) +x_{2}(s)) +\frac{1}{2}l_{xx}(s,\bar{x}(s),\bar{u}(s))(x_{1}(s) x_{1}(s))
\mathrm{d}s\bigg]  \notag \\
&&+\mathbb{E}\bigg[\int_{0}^{T}l(s,\bar{x}(s),u^{\varepsilon }(s)) -l(s,\bar{x}(s),\bar{u}(s))\mathrm{d}s\bigg]  \notag \\
&&+\mathbb{E}\bigg[h_{x}(\bar{x}(T))(x_{1}(T) +x_{2}(T)) +\frac{1}{2}h(\bar{x}(T),\alpha(T))(x_{1}(T)x_{1}(T))\bigg]  \geq  o(\varepsilon).  \label{varin}
\end{eqnarray}
\end{lemma}

\noindent
\textbf{Proof.}
From the fact that $(\bar{x}(\cdot),\bar{u}(\cdot),\alpha(\cdot))$ is optimal triple, we have
\begin{eqnarray*}
&&\mathbb{E}\bigg[\int_{0}^{T}l(s,x^{\varepsilon }(s),u^{\varepsilon }(s),\alpha(s)) \mathrm{d}s+h(x^{\varepsilon }(T),\alpha(T))\bigg] \\
&& \quad -\mathbb{E}\bigg[\int_{0}^{T}l(s,\bar{x}(s),\bar{u}(s),\alpha(s))\mathrm{d}s+h(\bar{x}(T),\alpha(T))\bigg] \geq 0.
\end{eqnarray*}
By Lemma \ref{est}, we get
\begin{eqnarray*}
0 &\leq &\mathbb{E}\bigg[\int_{0}^{T}\big[l(s,\bar{x}(s)+x_{1}(t) +x_{2}(t),u^{\varepsilon}(s),\alpha(s)) -l(s) \big ]\mathrm{d}s \\
&&+h(\bar{x}(T) +x_{1}(T) +x_{2}(T),\alpha(T)) -h(\bar{x}(T),\alpha(T)) \bigg]+o(\varepsilon) \\
&=&\mathbb{E}\bigg[\int_{0}^{T}\big[l(s,\bar{x}(s)+x_{1}(t) +x_{2}(t),\bar{u}(s),\alpha(s)) -l(s)\big] \mathrm{d}s\bigg] \\
&& \quad +\mathbb{E}\bigg [\int_{0}^{T}\big[l(s,\bar{x}(s)+x_{1}(t) +x_{2}(t),u^{\varepsilon }(s),\alpha(s)) \\
&& \quad -l(s,\bar{x}(s) +x_{1}(t) +x_{2}(t),\bar{u}(s),\alpha(s))\big]\mathrm{d}s\bigg] \\
&& \quad +\mathbb{E}\big[h(\bar{x}(T) +x_{1}(T)+x_{2}(T),\alpha(T))  -h(\bar{x}(T),\alpha(T))\big] +o(\varepsilon) \\
%\end{eqnarray*}
%%
%\begin{eqnarray*}
&=&\mathbb{E}\bigg [\int_{0}^{T}\big [l_{x}(s,\bar{x}(s),\bar{u}(s),\alpha(s))(x_{1}(s) +x_{2}(s)) \\
&&+\frac{1}{2}l_{xx}(s,\bar{x}(s),\bar{u}(s),\alpha(s))(x_{1}(s) +x_{2}(s))(x_{1}(s) +x_{2}(s))\big]\mathrm{d}s\bigg] \\
&&+\mathbb{E}\bigg[\int_{0}^{T}l(s,\bar{x}(s),u^{\varepsilon }(s),\alpha(s)) -l(s,\bar{x}(s),\bar{u}(s),\alpha(s))\mathrm{d}s\bigg] \\
&&+\mathbb{E}\bigg[\int_{0}^{T}(l_{x}(s,\bar{x}(s),u^{\varepsilon }(s),\alpha(s))-l_{x}(s,\bar{x}(s),\bar{u}(s),\alpha(s)))(x_{1}(s) +x_{2}(s)) \mathrm{d}s\bigg] \\
&&+\frac{1}{2}\mathbb{E}\bigg[\int_{0}^{T}(l_{xx}(s,\bar{x}(s),u^{\varepsilon }(s),\alpha(s)) -l_{xx}(s,\bar{x}(s),\bar{u}(s),\alpha(s))) x_{1}(s) x_{1}(s)\mathrm{d}s\bigg] \\
&&+\mathbb{E}\bigg[h_{x}(\bar{x}(T))(x_{1}(T) +x_{2}(T)) +\frac{1}{2}h_{xx}(\bar{x}(T)) x_{1}(T) x_{1}(T)\bigg] + o(\varepsilon).
\end{eqnarray*}
By Lemma \ref{est}, we complete the proof.
\hfill $\Box $

Next, we shall find the first and second-order adjoint equations. Using these processes, we are able to derive the variational inequality
from(\ref{varin}).

Consider the following stochastic system:
\begin{equation}
\left\{\begin{array}{rcl}
\mathrm{d}z(t) & = &(b_{x}(t)z(t) +\phi(t)) \mathrm{d}t+(\sigma_{x}(t) z t) +\psi(t))\mathrm{d}W(t) \\
&  & +(\gamma_{x}(t)z(t) +\mu(t))\mathrm{d}\tilde{V}(t) \\
z(0) & = & 0,
\end{array}\right.  \label{f1}
\end{equation}
where $(\phi(\cdot),\psi(\cdot),\mu(\cdot)) \in \mathcal{M}^{2}(0,T;\mathbb{R}^{L})\times\mathcal{M}^{2}(0,T;\mathbb{R}^{L\times d})\times \mathcal{V}^{2}(0,T;\mathbb{R}^{L\times D})$.
We construct a linear functional on the
Hilbert space $\mathcal{M}^{2}(0,T;\mathbb{R}^{L})\times \mathcal{M}^{2}(0,T;\mathbb{R}^{L\times d})\times \mathcal{V}^{2}(0,T;\mathbb{R}^{L\times D})$ as follows:
\begin{equation*}
\mathcal{I}_{1}(\phi(\cdot),\psi(\cdot),\mu(\cdot)) =\mathbb{E}\bigg[ \int_{0}^{T}l_{x}(t,\bar{x}(t),\bar{u}(t))z(t)+h_{x}(\bar{x}(T),\alpha(T)) z(T)\bigg],
\end{equation*}
where $\phi(\cdot),\psi(\cdot),\mu(\cdot)$ and $z(\cdot)$ are defined in(\ref{f1}).
Then, by virtue of the Riesz Representation Theorem, there exists a unique solution
\begin{equation*}
(p(\cdot),q(\cdot),s(\cdot)) \in \mathcal{M}^{2}(0,T;\mathbb{R}^{L})\times \mathcal{M}^{2}(0,T;\mathbb{R}^{L\times d})\times \mathcal{V}^{2}(0,T;\mathbb{R}^{L\times D})
\end{equation*}
such that
\begin{eqnarray*}
&&\mathbb{E}\bigg[\int_{0}^{T}(\langle p(t),\phi(t) \rangle +\langle q(t),\psi(t) \rangle +\sum_{m=1}^{D}\sum_{n=1}^{L}s_{nm}(t)
\mu_{nm}(t) \lambda_{m}(t)) \mathrm{d}t\bigg]
= \mathcal{I}_{1}(\phi(\cdot),\psi(\cdot),\mu(\cdot)),
\end{eqnarray*}
for $\forall(\phi(\cdot),\psi(\cdot),\mu(\cdot)) \in
\mathcal{M}^{2}(0,T;\mathbb{R}^{L})\times \mathcal{M}^{2}(0,T;
\mathbb{R}^{L\times d})\times \mathcal{V}^{2}(0,T;\mathbb{R}^{L\times D})$.
We now consider equations (\ref{v1}) and (\ref{v2}) as the solution to (\ref{f1}), respectively. Then, we have
\begin{eqnarray*}
&&\mathbb{E}\bigg[ \int_{0}^{T}l_{x}(t,\bar{x}(t),\bar{u}(t)) x_{1}(t) +h_{x}(\bar{x}(T),\alpha(T)) x_{1}(T)\bigg] \\
&=&\mathbb{E}\bigg[\int_{0}^{T}\Big(\langle p(t),\delta b^{\varepsilon }(t) \rangle +\langle q(t),\delta \sigma^{\varepsilon }(t) \rangle +\sum_{m=1}^{D}\sum_{n=1}^{L}s_{nm}(t) \cdot \delta \gamma_{nm}(t) \cdot \lambda_{m}(t) \Big)\mathrm{d}t\bigg],
\end{eqnarray*}
and
\begin{eqnarray*}
&&\mathbb{E}\bigg[ \int_{0}^{T}l_{x}(t,\bar{x}(t),\bar{u}(t)) x_{2}(t) +h_{x}(\bar{x}(T),\alpha(T)) x_{2}(T)\bigg] \\
&=&\mathbb{E}\bigg [\int_{0}^{T}\bigg(\langle p(t),\frac{1}{2}b_{xx}(t) x_{1}(t) x_{1}(t) +\delta b_{x}^{\varepsilon }(t) x_{1}(t) \rangle
+\langle q(t),\frac{1}{2}\sigma_{xx}(t)x_{1}(t) x_{1}(t) +\delta \sigma_{x}^{\varepsilon}(t) x_{1}(t) \rangle \\
&&+\langle s(t),\frac{1}{2}\gamma_{xx}(t)x_{1}(t) x_{1}(t) +\delta \gamma_{x}^{\varepsilon}(t) x_{1}(t) \rangle \bigg)\mathrm{d}t\bigg] \\
&=&\mathbb{E}\bigg [\int_{0}^{T}\Big(\frac{1}{2}(p(t)b_{xx}(t) +q(t),\sigma_{xx}(t)+s(t) \gamma_{xx}(t)) x_{1}(t)x_{1}(t) \\
&&+(p(t) \delta b_{x}^{\varepsilon}(t)+q(t) \delta \sigma_{x}^{\varepsilon}(t) +s(t)\delta\gamma_{x}^{\varepsilon }(t))x_{1}(t) \Big)\mathrm{d}t\bigg].
\end{eqnarray*}
Then, using Assumption (A1) and Lemma \ref{est}, we can re-write the expression (\ref{varin}) as
\begin{eqnarray}
&&\mathbb{E}\bigg[\int_{0}^{T}H(t,\bar{x}(t),u^{\varepsilon }(t),\alpha(t),p(t),q(t),s(t)) -H(t,\bar{x}(t),\bar{u}(t),\alpha(t),p(t),q(t),s(t))\mathrm{d}t\bigg]  \notag \\
&& \quad +\frac{1}{2}\mathbb{E}\bigg[ \int_{0}^{T}x_{1}^{\top }(t)H_{xx}(t,\bar{x}(t),\bar{u}(t),\alpha(t),p(t),q(t),s(t))x_{1}(t) \mathrm{d}t\bigg]  \notag \\
&& \quad +\frac{1}{2}\mathbb{E}\big[ x_{1}^{\top }(T) h_{xx}(\bar{x}(T),\alpha(T))(x_{1}(T))\big]  \geq o(\varepsilon).  \label{var2}
\end{eqnarray}
Next, we focus on the quadratic terms in (\ref{var2}) by employing the Riesz Representation Theorem. Applying It\^{o}'s formula to
$X(t) =x_{1}(t) x_{1}^{\top }(t)$
% \%begin{equation*} % =[
%\begin{array}{ccc}
%x_{1}^{1}(t) x_{1}^{1}(t) & \cdots & x_{1}^{1}(
%t) x_{1}^{n}(t) \\
%\vdots &  & \vdots \\
%x_{1}^{n}(t) x_{1}^{1}(t) & \cdots & x_{1}^{n}(
%t) x_{1}^{n}(t)%
%\end{array}%
%]
%\end{equation*}%
on $[ 0,T]$ yields
\begin{equation}
\left\{\begin{array}{rcl}
\mathrm{d}X(t) & = & \Big \{X(t) b_{x}^{\top}(t) +b_{x}(t) X(t) +\delta b^{\varepsilon}(t) x_{1}^{\top }(t) +x_{1}(t) \delta b^{\varepsilon }(t)^{\top } \\
&  & +\displaystyle\sum_{j=1}^{d}\Big[\sigma_{x}^{j}(t) X(t)\sigma_{x}^{j}(t)^{\top }+\sigma_{x}^{j}(t)x_{1}(t) \delta \sigma^{\varepsilon j}(t)^{\top} \\
&  & +\delta \sigma^{\varepsilon j}(t) x_{1}(t)\sigma_{x}^{j}(t)^{\top }+\delta \sigma^{\varepsilon j}(t) \delta \sigma^{\varepsilon j}(t)^{\top }\Big ] \\
&  & +\displaystyle\sum_{j=1}^{D}\Big[\gamma_{x}^{j}(t)X(t)\gamma_{x}^{j}(t)^{\top}+\gamma_{x}^{j}(t)x_{1}(t)\delta\gamma^{\varepsilon j}(t)^{\top} \\
&  & +\delta \gamma^{\varepsilon j}(t) x_{1}(t)\gamma_{x}^{j}(t)^{\top }+\delta \gamma^{\varepsilon j}(t) \delta \gamma^{\varepsilon j}(t)^{\top }\Big ]\lambda_{j}(t) \Big \}\mathrm{d}t \\
&  & +\displaystyle\sum_{j=1}^{d}\Big[\sigma_{x}^{j}(t) X(t)+X(t) \sigma_{x}^{j}(t)^{\top }  +x_{1}(t) \delta \sigma^{\varepsilon j}(t)^{\top }+\delta \sigma^{\varepsilon j}(t) x_{1}(t)^{\top }\Big)\mathrm{d}W^{j}(t) \\
&  & +\displaystyle\sum_{j=1}^{D}\Big[\gamma_{x}^{j}(t) X(t)+X(t) \gamma_{x}^{j}(t)^{\top }+\delta \gamma^{\varepsilon j}(t) x_{1}(t)^{\top }+x_{1}(t) \delta \gamma^{\varepsilon j}(t)^{\top} \\
&  & +\gamma_{x}^{j}(t) X(t) \gamma_{x}^{j}(t)^{\top }+\gamma_{x}^{j}(t) x_{1}(t) \delta\gamma^{\varepsilon j}(t)^{\top }+\delta\gamma^{\varepsilon
j}(t) x_{1}(t)^{\top }\gamma_{x}^{j}(t)^{\top} \\
&  & +\delta \gamma^{\varepsilon j}(t) \delta \gamma^{\varepsilon j}(t)^{\top}\Big]\mathrm{d}\tilde{\Phi}^{j}(t) \\
X(t) & = & 0.
\end{array}\right.  \label{qudra}
\end{equation}
Consider the following symmetric matrix-valued linear stochastic differential equations:
\begin{equation}
\left\{\begin{array}{rcl}
\mathrm{d}Z(t) & = & \Big[Z(t)b_{x}^{\top}(t) +b_{x}(t) Z(t) +\sum_{j=1}^{d}\sigma_{x}^{j}(t) Z(t) \sigma_{x}^{j}(t)^{\top} \\
&  & +\displaystyle\sum_{j=1}^{D}\gamma_{x}^{j}(t) Z(t) \gamma_{x}^{j}(t)^{\top }\lambda_{j}(t) +\Upsilon(t) \Big ]\mathrm{d}t \\
&  & +\displaystyle\sum_{j=1}^{d}(\sigma_{x}^{j}(t) Z(t)+Z(t) \sigma_{x}^{j}(t)^{\top }+\Psi^{j}(t)) \mathrm{d}W^{j}(t) \\
&  & +\displaystyle\sum_{j=1}^{D}(\gamma_{x}^{j}(t) Z(t)+Z(t) \gamma_{x}^{j}(t)^{\top }+\gamma_{x}^{j}(t) Z(t) \gamma_{x}^{j}(t)^{\top }+\Pi^{j}(t)) \mathrm{d}\tilde{\Phi}^{j}(t) \\
Z(0) & = & 0,
\end{array}\right.  \label{f2}
\end{equation}
where $(\Upsilon(\cdot),\Psi(\cdot),\Pi(\cdot)) \in \mathcal{M}^{2}(0,T;\mathbb{R}^{L\times L})\times(\mathcal{M}^{2}(0,T;\mathbb{R}^{L\times L}))^{d}\times(\mathcal{V}^{2}(0,T;\mathbb{R}^{L\times L}))^{D}$.

Define a linear functional based on (\ref{f2}):
\begin{equation*}
\mathcal{I}_{2}(\Upsilon(\cdot),\Psi(\cdot),\Pi(\cdot)) =\mathbb{E}\bigg[\int_{0}^{T}\langle Z(t),H_{xx}(t)\rangle \mathrm{d}t+\langle h_{xx}(\bar{x}(T),\alpha(T)),Z(T) \rangle\bigg].
\end{equation*}
Obviously, $\mathcal{I}_{2}(\cdot,\cdot,\cdot)$
is a linear continuous functional on $\mathcal{M}^{2}(0,T;\mathbb{R}^{L\times L})\times(\mathcal{M}^{2}(0,T;\mathbb{R}^{L\times L}))^{d}$ $\times(\mathcal{V}^{2}(0,T;\mathbb{R}^{L\times L}))^{D}$.
There exists a unique solution $(P(\cdot),Q(\cdot),S(\cdot)) \in \mathcal{M}^{2}(0,T;\mathbb{R}^{L\times L})\times(\mathcal{M}^{2}(0,T;\mathbb{R}^{L\times L}))^{d}\times
(\mathcal{V}^{2}(0,T;\mathbb{R}^{L\times L}))^{D}$, such that
\begin{eqnarray}
\mathcal{I}_{2}(\Upsilon(\cdot),\Psi(\cdot),\Pi(\cdot)) &=&\mathbb{E}\bigg [\int_{0}^{T}\Big(\langle P(t),\Upsilon(\cdot)\rangle +\sum_{j=1}^{d}\langle Q^{j}(t),\Psi^{j}(t) \rangle  \notag \\
&& \quad +\sum_{j=1}^{D}\sum_{m=1}^{L}\sum_{n=1}^{L}S_{nm}^{j}(t) \Pi_{nm}^{j}(t) \lambda_{m}(t) \Big)\mathrm{d}t\bigg].  \label{forsec}
\end{eqnarray}
Now let
\begin{eqnarray*}
\Upsilon^{\varepsilon }(\cdot) &=&\delta b^{\varepsilon}(t) x_{1}^{\top }(t) +x_{1}(t) \delta b^{\varepsilon }(t)^{\top }
+\sum_{j=1}^{d}\Big[\sigma_{x}^{j}(t) x_{1}(t) \delta \sigma^{\varepsilon j}(t)^{\top } \\
&&+\delta \sigma^{\varepsilon j}(t) x_{1}(t) \sigma_{x}^{j}(t)^{\top }+\delta \sigma^{\varepsilon j}(t) \delta \sigma^{\varepsilon j}(t)^{\top }\Big ] \\
&&+\sum_{j=1}^{D}\Big [\gamma_{x}^{j}(t) x_{1}(t)\delta \gamma^{\varepsilon j}(t)^{\top }+\delta \gamma^{\varepsilon j}(t) x_{1}(t) \gamma_{x}^{j}(t)^{\top } +\delta \gamma^{\varepsilon j}(t) \delta \gamma^{\varepsilon j}(t)^{\top }\Big ]\lambda_{j}(t), \\
\Psi^{\varepsilon j}(t) &=&\sum_{j=1}^{d}x_{1}(t)\delta \sigma^{\varepsilon j}(t)^{\top }+\delta \sigma^{\varepsilon j}(t) x_{1}(t), \\
\Pi^{\varepsilon j}(t) &=&\sum_{j=1}^{D}\delta \gamma^{\varepsilon j}(t) x_{1}(t)^{\top }+x_{1}(t) \delta \gamma^{\varepsilon j}(t)^{\top }+\gamma_{x}^{j}(t) x_{1}(t) \delta \gamma^{\varepsilon j}(t)^{\top} \\
&&+\delta \gamma^{\varepsilon j}(t) x_{1}(t)^{\top}\gamma_{x}^{j}(t)^{\top }+\delta \gamma^{\varepsilon j}(t) \delta \gamma^{\varepsilon j}(t)^{\top}.
\end{eqnarray*}
Then, the relation (\ref{var2}) can be expressed as
\begin{eqnarray}
&&\mathbb{E}\bigg[\int_{0}^{T}H(t,\bar{x}(t-),u^{\varepsilon }(t),\alpha(t-),p(t-),q(t),s(t)) \notag \\
&& \quad -H(t,\bar{x}(t-),\bar{u}(t),\alpha(t-),p(t-),q(t),s(t))\mathrm{d}t\bigg] \notag \\
&& \quad +\frac{1}{2}\mathbb{E}\bigg [\int_{0}^{T}\Big(\langle P(t),\Upsilon^{\varepsilon }(\cdot) \rangle
+\sum_{j=1}^{d}\langle Q^{j}(t),\Psi^{\varepsilon j}(t) \rangle \mathrm{d}t  \notag \\
&& \quad +\sum_{j=1}^{D}\sum_{m=1}^{L}\sum_{n=1}^{L}S_{nm}^{j}(t) \Pi_{nm}^{\varepsilon j}(t) \lambda_{m}(t) \Big)\mathrm{d}t\bigg ]  \geq o(\varepsilon),  \label{var3}
\end{eqnarray}

From the definition of $(\Upsilon^{\varepsilon }(\cdot),\Psi^{\varepsilon j}(\cdot),\Pi^{\varepsilon j}(\cdot)) $ and Lemma \ref{est}, we have
\begin{eqnarray}
&&\mathbb{E}\bigg [\int_{0}^{T}H(t,\bar{x}(t-),u^{\varepsilon }(t),\alpha(t-),p(t-),q(t),s(t))  \notag \\
&& \quad -H(t,\bar{x}(t-),\bar{u}(t),\alpha(t-),p(t-),q(t),s(t))\mathrm{d}t\bigg ]  \notag \\
&& \quad +\frac{1}{2}\mathbb{E}\bigg [\int_{0}^{T}\Big(\text{tr}[ \delta\sigma^{\varepsilon }(t)^{\top }P(t) \cdot \delta
\sigma^{\varepsilon }(t) ] +\sum_{j=1}^{D}\delta \gamma^{\varepsilon j}(t)^{\top }P(t) \delta \gamma^{\varepsilon j}(t) \lambda_{j}(t)  \notag \\
&& \quad +\sum_{j=1}^{D}\sum_{m=1}^{L}\sum_{n=1}^{L}(\delta \gamma^{\varepsilon }(t)^{\top }S^{j}(t) \delta \gamma^{\varepsilon }(t))_{mn}\lambda_{j}(t) \Big)\mathrm{d}t\bigg]  \geq o(\varepsilon),  \label{var4}
\end{eqnarray}
where
\begin{eqnarray*}
\delta \sigma^{\varepsilon }(t) &=&\sigma(t,\bar{x}(t),u^{\varepsilon}(t),\alpha(t)) -\sigma(t,\bar{x}(t),\bar{u}(t),\alpha(t)), \\
\delta \gamma^{\varepsilon j}(t) &=&\gamma^{j}(t,\bar{x}(t),u^{\varepsilon}(t),\alpha(t)) -\gamma^{j}(t,\bar{x}(t),\bar{u}(t),\alpha(t)).
\end{eqnarray*}
Moreover, we adopt the scalar product: $\langle
A,B\rangle =$tr$[ AB]$, $A$, $B\in \mathbb{R}^{n\times n}$. It therefore follows that for all $y\in \mathbb{R}^{n}$ and $A\in \mathbb{R}^{n\times n}$, we have $\langle yy^{\top },A\rangle =$tr$[(yy^{\top }) A] =y^{\top }Ay$.

\vspace{0.25cm}

\begin{remark}
The inequality (\ref{var4}) is called the variational inequality of our
optimal control problem. In contrast to the classical work by Peng \cite%
{Peng1990}, the appearance of extra term with $\lambda_{j}(t) $
is due to the random environments.
\end{remark}

Next we shall find the second-order adjoint equation. The first one has been
given in Zhang et al \cite{ZEK} like:
\begin{equation}
\left\{\begin{array}{rcl}
-\mathrm{d}p(t) & = & H_{x}(t,\bar{x}(t-),\bar{u}(t),\alpha(t-),p(t-),q(t),s(t)) \mathrm{d}t \\
&  & -q(t) \mathrm{d}W(t) -s(t) \mathrm{d}\tilde{\Phi}(t), \\
p(T) & = & h_{x}(\bar{x}(T),\alpha(T)).
\end{array}\right.  \label{adj1}
\end{equation}
From (\ref{forsec}) and (\ref{var4}), we can find the second-order adjoint equation is
\begin{equation}
\left\{\begin{array}{rcl}
-\mathrm{d}P(t) & = & \bigg[P(t) b_{x}(t)+b_{x}(t)^{\top }P(t) \mathrm{d}t+\displaystyle\sum_{j=1}^{d}\sigma_{x}^{j\top }(t) P(t) \sigma_{x}^{j}(t) \\
&  & +\displaystyle\sum_{j=1}^{d}(\sigma_{x}^{j\top }(t) Q^{j}(t) +Q^{j}(t) \sigma_{x}^{j}(t)) \\
&  & +\displaystyle\sum_{j=1}^{D}\gamma_{x}^{j\top }(t) S^{j}(t)+S^{j}(t) \gamma_{x}^{j}(t) +\gamma_{x}^{j}(t)^{\top }[ P(t) +S^{j}(t)]\gamma_{x}^{j}(t) \\
&  & +H_{xx}(t,\bar{x}(t-),\bar{u}(t),\alpha(t-),p(t-),q(t),s(t))\bigg]\mathrm{d}t \\
&  & -Q(t) \mathrm{d}W(t) -S(t) \mathrm{d}\tilde{\Phi}(t), \\
P(T) & = & h_{xx}(\bar{x}(T),\alpha(T)),
\end{array}\right.  \label{adj2}
\end{equation}
where $H$ is defined in (\ref{Hal}). Now define an $\mathcal{H}$-function:
\begin{eqnarray*}
\mathcal{H}(t,x,u,e_{i}) & \triangleq & H(t,x,u,e_{i},p(t),q(t),s(t)) \\
&&-\frac{1}{2}\text{tr}[ \sigma(t,\bar{x}(t),\bar{u}(t),e_{i})^{\top }P(t) \sigma(t,\bar{x}(t),\bar{u}(t),e_{i})] \\
&&+\frac{1}{2}\text{tr}\big\{[\delta \sigma^{\varepsilon}(t)]^{\top }P(t)[\delta \sigma^{\varepsilon}(t)]\big\} \\
&&-\frac{1}{2}\sum_{j=1}^{D}\gamma^{j}(t,\bar{x}(t),\bar{u}(t),e_{i})^{\top }P(t) \gamma^{j}(t,\bar{x}(t),\bar{u}(t),e_{i}) \lambda_{j}(t) \\
&&+\frac{1}{2}\sum_{j=1}^{D}\text{tr}\big\{[\delta \gamma^{\varepsilon }(t)]^{\top }P(t)[\delta \gamma^{\varepsilon}(t)]\lambda_{j}(t)\big\} \\
&&-\frac{1}{2}\sum_{j=1}^{D}\sum_{m=1}^{L}\sum_{n=1}^{L}\Big(\delta \gamma^{\varepsilon }(t)^{\top}S^{j}(t) \delta\gamma^{\varepsilon}(t) \\
&&-\gamma^{j}(t,\bar{x}(t),\bar{u}(t),e_{i})^{\top }S^{j}(t) \gamma^{j}(t,\bar{x}(t),\bar{u}(t),e_{i})\Big)\lambda_{j}(t).
\end{eqnarray*}
Inequality (\ref{var4}) can be expressed as
\begin{equation}
\mathcal{H}(t,\bar{x}(t),\bar{u}(t),e_{i}) =\inf_{u\in U}\mathcal{H}(t,\bar{x}(t),u,e_{i}).  \label{var5}
\end{equation}

\begin{theorem}\label{the} \sl
Suppose that Assumptions \emph{(A1)-(A3)} are in force. If $(\bar{x}(\cdot),\bar{u}(\cdot)) $ is a
solution to optimal control problem(\ref{cpr}). Then
\begin{eqnarray*}
(p(\cdot),q(\cdot),s(\cdot)) & \in & \mathcal{M}^{2}(0,T;\mathbb{R}^{L})\times \mathcal{M}^{2}(0,T;\mathbb{R}^{L\times d})\times \mathcal{V}^{2}(0,T;\mathbb{R}^{L\times D}) \\
(P(\cdot),Q(\cdot),S(\cdot)) & \in & \mathcal{M}^{2}(0,T;\mathbb{R}^{L\times L})\times(\mathcal{M}^{2}(0,T;\mathbb{R}^{L\times L}))^{d}\times(\mathcal{V}^{2}(0,T;\mathbb{R}^{L\times L}))^{D},
\end{eqnarray*}
which are, respectively, unique solutions to equations \emph{(\ref{adj1})-(\ref{adj2})} such that the variational inequality \emph{(\ref{var4})} holds.
\end{theorem}

\noindent
\textbf{Proof.}
From (\ref{var4}) and discussion above, we complete the proof.
\hfill $\Box$

We now establish the relationship between MP and DPP. To this end, we consider the Markovian (feedback) control, namely, the control $u(t)$ of the form $u(t,X(t-),\alpha(t-))$.

Define
\begin{equation*}
J(t,x,e_{i},u) =\mathbb{E}\bigg[\int_{t}^{T}l(s,x(s),u(s),\alpha(s)) \mathrm{d}s+h(x(T),\alpha(T))\bigg], \quad \text{for } u \in \mathcal{U}_{ad}.
\end{equation*}
Let us define the \textit{value function} as follows:
\begin{equation}
V(t,x,e_{i}) =\inf_{u\in \mathcal{U}_{ad}}J(
t,x,e_{i},u)  \label{optpro}
\end{equation}
for every $(t,x,e_{i}) \in [ 0,T] \times \mathbb{R}^{L}\times S$.

Employing the classical dynamic programming principle (see Fleming and Soner \cite{FS}), we are able to attain the following Hamilton-Jacobi-Bellman (HJB, for short) equation with terminal boundary condition:
\begin{equation}
\left\{\begin{array}{l}
\displaystyle\frac{\partial }{\partial t}V(t,x,e_{i}) +\inf_{u\in U}\{\mathcal{L}V(t,x,e_{i}) +l(t,x,u,e_{i}) \} =0 \\
V(T,x,e_{i}) =h(x,e_{i}),
\end{array}\right.  \label{HJB}
\end{equation}
where
\begin{eqnarray*}
\mathcal{L}^{u}V(t,x,e_{i}) &=&\sum_{j=1}^{L}\frac{\partial V(t,x,e_{i}) }{\partial x_{j}}b_{i}(t,x,u,e_{i}) \\
&&+\frac{1}{2}\sum_{i=1}^{L}\sum_{j=1}^{L}\frac{\partial^{2}V(t,x,e_{i}) }{\partial x_{i}\partial x_{j}}\sum_{l=1}^{d}(\sigma_{il}\sigma_{il})(t,x,u,e_{i}) \\
&&+\sum_{m=1}^{D}\Big [V(t,x+\gamma^{m}(t,x,u,e_{i}),e_{m}) -V(t,x,e_{i}) \\
&&-\sum_{j=1}^{L}\frac{\partial V(t,x,e_{i}) }{\partial x_{j}}\gamma_{jm}(t,x,u,e_{i}) \Big ]\lambda_{im}.
\end{eqnarray*}
Put
\begin{equation*}
G(t,x,u,p,P,e_{i}) =\mathcal{L}^{u}V(t,x,e_{i})-l(t,x,u,e_{i})
\end{equation*}
for $(t,x,p,P,e_{i}) \in [ 0,T] \times\mathbb{R}^{L}\times U\times \mathbb{R}^{L}\times \mathbb{R}^{L\times L}\times S$.

\vspace{0.25cm}

\begin{remark}
Clearly, because of the presence of regime-switching, (\ref{HJB}) is a
coupled system of nonlinear differential equations. Meanwhile, it is
impossible for us to obtain the explicit form of the value function and an
optimal control by solving (\ref{HJB}). The main reason comes from the fact
that (\ref{HJB}) is a coupled system of nonlinear second-order partial
differential equations, adding extreme difficulty in finding a closed-form
solution of (\ref{HJB}). Recall in lieu of a single nonlinear differential
equation in the traditional literature for optimal control problems, we need
more powerful tools and analysis to deal with it.
\end{remark}

As a matter of fact, one should seek an optimal control and associated state
trajectories. Normally, the main motivation for studying the dynamic
principle programming is that one might construct an optimal control through
the value function. The following result gives a way of testing whether a
given admissible pair is optimal and, more importantly, suggests how to
construct an optimal feedback control. Such a result is called a \textit{%
verification theorem}. For smooth case, we have a following result.

\begin{theorem}[Verification theorem]\label{verfica} \sl
Assume that \emph{(A1)-(A3)} hold. Let $V\in C^{1,2}([0,T] \times \mathbb{R}^{n}\times e_{i})$ for each $e_{i}\in S$ be a
solution of HJB equation (\ref{HJB}). Then
\begin{equation}
V(t,x,e_{i}) \leq J(t,x,e_{i}), \quad \forall u(\cdot) \in \mathcal{U}_{ad}[t,T],  \quad (t,x) \in [ 0,T] \times \mathbb{R}^{n}.  \label{verf1}
\end{equation}
Furthermore, an admissible pair $(\bar{x}(\cdot),\bar{u}(\cdot)) $ is an optimal control for(\ref{optpro}) if and only if
\begin{eqnarray}
&&\frac{\partial }{\partial t}V(s,\bar{x}(s),e_{i})+\inf_{u\in U}\{ \mathcal{L}^{u}V(s,\bar{x}(s),e_{i}) +l(s,\bar{x}(s),u,e_{i})\} \notag \\
&=&\frac{\partial }{\partial t}V(s,\bar{x}(s),e_{i}) +\mathcal{L}^{\bar{u}}V(s,\bar{x}(s),e_{i}) +l(s,\bar{x}(s),\bar{u}(s),e_{i})  \notag \\
&=&0,  \quad t\leq s\leq T.  \label{verf2}
\end{eqnarray}
\end{theorem}

\noindent
\textbf{Proof.}
Fix $e_{i}\in S$. Let $u(\cdot) \in \mathcal{U}_{ad}[t,T]$ with its associated state trajectory $x(\cdot)$.
Applying It\^{o}' formula to $V(t,x,e_{i}) $ on $[ t,T]$, we have
\begin{eqnarray}
V(t,x,e_{i}) &=&\mathbb{E}\bigg[h(x(T),\alpha(T)) +\int_{t}^{T}l(s,x(s),u(s),\alpha(s)) \mathrm{d}s\bigg]  \notag \\
&&-\mathbb{E}\bigg[\int_{t}^{T}V_{t}(s,x(s),e_{i}) +G(s,x(s),V_{x}(s,x(s),e_{i}),V_{xx}(s,x(s),e_{i}),e_{i}) \mathrm{d}s\bigg]  \notag \\
&\leq &J(t,x,e_{i}) +\mathbb{E}\bigg [\int_{t}^{T}-V_{t}(s,x(s),e_{i})  \notag \\
&&+G(s,x(s),u(s),-V_{x}(s,x(s),e_{i}),-V_{xx}(s,x(s),e_{i}),e_{i}) \mathrm{d}s\bigg ]  \notag \\
&\leq &J(t,x,e_{i}) +\mathbb{E}\bigg [\int_{t}^{T}-V_{t}(s,x(s),e_{i})  \notag \\
&&+\sup_{u\in U}G(s,x(s),u,-V_{x}(s,x(s),e_{i}),-V_{xx}(s,x(s),e_{i}),e_{i}) \mathrm{d}s\bigg ]  \notag \\
&=&J(t,x,e_{i}),  \label{verf3}
\end{eqnarray}
which implies (\ref{verf1}). Applying (\ref{verf3}) to $(\bar{x}(
\cdot),\bar{u}(\cdot))$, together with (\ref{verf2}), yields the desired result.
\hfill $\Box$

\vspace{0.25cm}

\begin{remark}
As claimed in Theorem \ref{verfica}, the value function requires to be
smooth enough. As we have well known, the value function \ref{optpro} is not
necessarily smooth. Thus, one perhaps introduce the viscosity solution to
study our control problems. As this complete remark of the existence is much
longer than the present paper, it will be reported elsewhere.
\end{remark}

\label{App}
Let us now provide a substantial example to demonstrate how to construct an optimal control by Theorem \ref{the}.

\begin{example}
Consider the following control system $(L=d=1)$:
\begin{equation}
\left\{\begin{array}{rcl}
\mathrm{d}x(t) & = & u(t) \sigma(\alpha(t-)) \mathrm{d}W(t) +u(t) \gamma^{\top }(\alpha(t-)) \mathrm{d}\tilde{\Phi}(t), \\
x(0-) & = & 0,
\end{array}\right.  \label{exa}
\end{equation}
where $\sigma :\mathbb{R}^{D} \rightarrow \mathbb{R}$, $\gamma :\mathbb{R}^{D} \rightarrow \mathbb{R}^{D}$ with the control domain being $[ 0,1]$ and the cost functional being
\begin{equation}
J(u(\cdot)) =\mathbb{E}\bigg[ -\int_{0}^{T}u(t) \mathrm{d}t+\frac{1}{2}x^{2}(T) +\alpha^{\top }(T) \alpha(T)\bigg].  \label{exacost}
\end{equation}
Let $(\bar{x}(\cdot),\bar{u}(\cdot))
$ be optimal pair to be determined. The first and second-order adjoint
equations are
\begin{equation}
\left\{\begin{array}{rcl}
\mathrm{d}p(t) & = & q(t) \mathrm{d}W(t)+\varsigma(t) \mathrm{d}\tilde{\Phi}(t), \\
p(T) & = & x(T),
\end{array}\right.  \label{exaadj1}
\end{equation}
and
\begin{equation}
\left\{\begin{array}{rcl}
\mathrm{d}P(t) & = & Q(t) \mathrm{d}W(t)+S(t) \mathrm{d}\tilde{\Phi}(t), \\
P(T) & = & 1.
\end{array}\right.  \label{exaadj2}
\end{equation}
According to the existence and uniqueness of BSDE (\ref{exaadj2}), the
unique adapted the solution is $(P(t),Q(t),S(t)) =(1,0,0)$, $t\in [ 0,T]$.
The corresponding $\mathcal{H}$-function is
\begin{eqnarray*}
\mathcal{H}(s,\bar{x}(s),u,e_{i}) &=&\frac{1}{2}u^{2}\sigma^{2}(e_{i}) - \big(1+\bar{u}(s) \sigma^{2}(e_{i}) -\sigma(e_{i}) q(s)\big)u \\
&&+\frac{1}{2}\sum_{j=1}^{D}u^{2}\gamma^{j}(\alpha(e_{i}))^{2}-u\sum_{j=1}^{D}\gamma^{j}(\alpha(e_{i}))^{2}\bar{u}(s) +u\sum_{j=1}^{D}\gamma^{j}(e_{i}) \varsigma_{j}\lambda_{im} \\
&=&\frac{1}{2}u^{2}\Big(\sigma^{2}(e_{i})+\sum_{j=1}^{D}\gamma^{j}(\alpha(e_{i}))^{2}\Big) \\
&&-\Big(1+\bar{u}(s) \sigma^{2}(e_{i}) -\sigma(e_{i}) q(s) +\bar{u}(s)\sum_{j=1}^{D}\gamma^{j}(\alpha(e_{i}))^{2} \\
&&-\sum_{j=1}^{D}\gamma^{j}(e_{i}) \varsigma_{j}\lambda_{im}\Big)u.
\end{eqnarray*}
The function $\mathcal{H}(s,\bar{x}(s),\cdot,e_{i}) $ attains its minimum at
\begin{equation*}
u=1+\bar{u}(s) \sigma^{2}(e_{i}) - \sigma(e_{i}) q(s) +\bar{u}(s) \sum_{j=1}^{D}\gamma^{j}(\alpha(e_{i}))^{2}-\sum_{j=1}^{D}\gamma^{j}(e_{i}) \varsigma_{j}\lambda_{ij}.
\end{equation*}
Particularly, set $\sigma(e_{i}) =1$, $\gamma(e_{i}) \equiv 0$. Then
\begin{equation*}
u=1+\bar{u}(s) -q(s) .
\end{equation*}
Hence the necessary condition of optimality as specified by the stochastic
maximum principle (Theorem \ref{the}) would be satisfied if one could find a
control $u(\cdot) $ such that the corresponding $q(s) =1$, $\varsigma(s) =0$. By the first-order adjoint
equation (\ref{exaadj1}), it is clear that if we take $u(s) =1$
with the corresponding state $\bar{x}(s) =W(s) $,
then the unique solution of (\ref{exaadj1}) is $(p(s),q(s),\varsigma(s)) =(W(s),1,0) $.
\end{example}

%In this section, we study a kind of linear quadratic portfolio optimization
%problem in a Markov regime-switching diffusion financial market. The optimal
%portfolio in the state feedback form can be obtained by applying the maximum
%principle(Theorem \ref{the}).
%
%Assume that the Markov regime-switching diffusion financial market involoves
%one risk-free asset and $L$ risky assets. The risk-free asset's price
%process $S_{0}(t) $ is given by:%
%\begin{equation*}
%\{
%\begin{array}{lll}
%\mathrm{d}S_{0}(t) & = & r(t,\alpha(t-)
%) S_{0}(t) \mathrm{d}t, \\
%S_{0}(0) & = & \chi >0,%
%\end{array}%
%.
%\end{equation*}%
%where $r(t,e_{i}) \geq 0$, $i=1,\ldots,D$, are bounded,
%deterministic functions on $[0,T]$ and can be interpreted as the interest
%rates in different market modes.
%
%The price of the other $L$ risky assets $S_{k}(t) $, $%
%k=1,2,\ldots,L$, are governed by the following forward stochastic
%differential equations:%
%\begin{equation*}
%\mathrm{d}S_{k}(t) =S_{k}(t-) [ b_{k}(
%t,\alpha(t-)) \mathrm{d}t+\sum_{j=1}^{d}\sigma
%_{kj}(t,\alpha(t-)) \mathrm{d}W^{j}(t) %
%] .
%\end{equation*}%
%Moreover, we shall suppose that $b_{k}(t,e_{i}) \geq $ $r(
%t,e_{i}) $

\section{Optimal control for recursive utilities}

Based on the results obtained in Section \ref{section3}, we study the
optimal control problem for stochastic recursive utilities for systems
composed of SDE (\ref{sde1}) and BSDE (\ref{bsde1}). The cost functional $J(u(\cdot))$ now is defined in (\ref{costfbsde}).
The control problem is to minimize $J(u(\cdot))$
over $\mathcal{U}_{ad}$. For simplicity, we assume that $L=1$, $d=1$. Let $\bar{u}(\cdot)$ be the optimal control and let $(\bar{x}(\cdot),\bar{y}(\cdot),\bar{z}(\cdot),\bar{\kappa}(\cdot))$ be the corresponding solution of the equations (\ref{sde1}) and (\ref{bsde1}). Similarly, we
define $(x^{\varepsilon }(\cdot),y^{\varepsilon}(\cdot),z^{\varepsilon }(\cdot),\kappa^{\varepsilon}(\cdot))$ for $u^{\varepsilon }(\cdot)$.

We shall give the variational equation for BSDE (\ref{bsde1}). For this aim,
we consider the following two adjoint equations:
\begin{equation}
\left\{\begin{array}{rcl}
-\mathrm{d}\mathfrak{p}(t) & = & F(t) dt-\mathfrak{q}(t) \mathrm{d}W(t) -\mathfrak{s}(t)\mathrm{d}\tilde{\Phi}(t), \\
\mathfrak{p}(t) & = & g_{x}(\bar{x}(T),\alpha(T))
\end{array}\right.  \label{fbsdeadj1}
\end{equation}
and
\begin{equation}
\left\{\begin{array}{rcl}
-\mathrm{d}\mathcal{P}(t) & = & G(t) dt-\mathcal{Q}(t) \mathrm{d}W(t) -\mathcal{S}(t)\mathrm{d}\tilde{\Phi}(t), \\
\mathcal{P}(T) & = & g_{xx}(\bar{x}(T),\alpha(T)),
\end{array}\right.  \label{fbsdeadj2}
\end{equation}
where $F(t)$ and $G(t)$ are adapted processes with suitable properties, will be determined later.

Applying It\^{o}'s formula to $\langle \mathfrak{p}(t),x_{1}(t) +x_{2}(t) \rangle +\frac{1}{2}x_{2}(t)^{\top }\mathcal{P}(t) x_{2}(t)$ yields
\begin{eqnarray*}
&&\langle \mathfrak{p}(T),x_{1}(T)+x_{2}(T) \rangle +\frac{1}{2}x_{2}(T)^{\top }\mathcal{P}(T) x_{2}(T) \\
&=&\langle \mathfrak{p}(t),x_{1}(t)+x_{2}(t) \rangle +\frac{1}{2}x_{2}(t)^{\top }\mathcal{P}(t) x_{2}(t) \\
&&+\int_{t}^{T}(A_{1}(s) +A_{2}(s)(x_{1}(s) +x_{2}(s)) +\frac{1}{2}A_{3}(s) x_{1}(s) x_{1}(s) +A_{4}(s)x_{1}(s)) \mathrm{d}s \\
&&+\int_{t}^{T}\Big [\mathfrak{p}(s) \delta \sigma^{\varepsilon}(s) +A_{5}(s)(x_{1}(s)+x_{2}(s)) \\
&&+\frac{1}{2}A_{6}(s) x_{1}(s) x_{1}(s)+A_{7}(s) x_{1}(s) \Big ]\mathrm{d}W(s)  \\
&&+\int_{t}^{T}\sum_{j=1}^{D}\Big [(\mathfrak{p}(s) +\mathfrak{s}_{j}(s)) \delta \gamma_{j}^{\varepsilon}(s) +\frac{1}{2}(\mathcal{P}(s) +\mathcal{S}(s)) \delta \gamma_{j}^{\varepsilon }(s)\delta \gamma_{j}^{\varepsilon }(s) \\
&&+A_{8}^{j}(s)(x_{1}(s) +x_{2}(s)) +\frac{1}{2}A_{9}^{j}(s) x_{1}(s)x_{1}(s) +A_{10}^{j}(s) x_{1}(s) \Big]\mathrm{d}\tilde{\Phi}_{j}(s),
\end{eqnarray*}
where
\begin{eqnarray*}
A_{1}(s) &=&\mathfrak{p}(s) \delta b^{\varepsilon}(s) +\delta \sigma^{\varepsilon }(s) \mathfrak{q}(s) +\frac{1}{2}\mathcal{P}(s) \delta \sigma^{\varepsilon }(s) \delta \sigma^{\varepsilon }(s) \\
&&+\sum_{j=1}^{D}(\delta \gamma_{j}^{\varepsilon }(s)\mathfrak{s}_{j}(s) +\frac{1}{2}(\mathcal{P}(s) +\mathcal{S}(s)) \delta \gamma_{j}^{\varepsilon}(s) \delta \gamma_{j}^{\varepsilon }(s))\lambda_{j}(s), \\
A_{2}(s) &=&\mathfrak{p}(s) b_{x}(s) +\mathfrak{q}(s) \sigma_{x}(s) +\sum_{j=1}^{D}\gamma_{x}^{j}(s) \mathfrak{s}_{j}(s) \lambda_{j}(s) -F(s), \\
A_{3}(s) &=&b_{xx}(s) \mathfrak{p}(s) +\mathfrak{q}(s) \sigma_{xx}(s) +2\mathcal{P}(
s) b_{x}(s) +\mathcal{P}(s) \sigma_{x}(s) \sigma_{x}(s) \\
&&+\mathcal{P}(s) \gamma_{x}(s) \gamma_{x}(s) \lambda_{j}(s) +2\mathcal{Q}(s) \sigma_{x}(s) -G(s) \\
&&+\sum_{j=1}^{D}(\gamma_{xx}^{j}(s) \mathfrak{s}_{j}(s) +\mathcal{S}(s) \gamma_{x}^{j}(s) \gamma_{x}^{j}(s) +2\mathcal{S}(s) \gamma_{x}^{j}(s)) \lambda_{j}(s), \\
A_{4}(s) &=&\mathfrak{p}(s) \delta b_{x}^{\varepsilon }(s) +\mathfrak{q}(s) \delta\sigma_{x}^{\varepsilon }(s) +\mathcal{P}(s) \delta b^{\varepsilon }(s) +\mathcal{Q}(s) \delta \sigma^{\varepsilon }(s) \\
&&+\sum_{j=1}^{D}\Big(\delta \gamma_{x}^{\varepsilon j}(s)\mathfrak{s}_{j}(s) +\mathcal{P}(s) \gamma_{x}^{j}(s) \delta \gamma_{mn}^{\varepsilon j}(s) \\
&&+\mathcal{S}(s) \delta \gamma_{j}^{\varepsilon }(s) +\mathcal{S}(s) \delta \gamma_{j}^{\varepsilon }(s) \gamma_{x}^{j}(s) \Big)\lambda_{j}(s) +\mathcal{P}(s) \delta \sigma^{\varepsilon }(s), \\
A_{5}(s) &=&\mathfrak{p}(s) \sigma_{x}(s) +\mathfrak{q}(s), \\
%\end{eqnarray*}
%%
%\begin{eqnarray*}
A_{6}(s) &=&\mathfrak{p}(s) \sigma_{xx}(s) +\mathcal{Q}(s) +2\mathcal{P}(s) \sigma_{x}(s), \\
A_{7}(s) &=&\mathfrak{p}(s) \delta \sigma_{x}^{\varepsilon }(s) +\mathcal{P}(s) \delta \sigma^{\varepsilon }(s), \\
A_{8}^{j}(s) &=&\mathfrak{p}(s) \gamma_{x}^{j}(s) +\mathfrak{s}_{j}(s) +\delta \gamma_{x}^{j}(s), \\
A_{9}^{j}(s) &=&\mathfrak{p}(s) \gamma_{xx}^{j}(s) +\mathfrak{s}_{j}(s) \gamma_{xx}^{j}(s) +2\mathcal{P}(s) \gamma_{x}^{j}(s) +\mathcal{P}(s) \gamma_{x}^{j}(s) \gamma_{x}^{j}(s) \\
&&+\mathcal{S}(s) +2\mathcal{S}(s) \gamma_{x}^{j}(s) +\mathcal{S}(s) \gamma_{x}^{j}(s) \gamma_{x}^{j}(s), \\
A_{10}^{j}(s) &=&(\mathfrak{p}(s) +\mathfrak{s}_{j}(s)) \delta \gamma_{x}^{\varepsilon j}(s) +\mathcal{P}(s) \delta \gamma_{j}^{\varepsilon }(s) +\mathcal{P}(s) \delta \gamma_{j}^{\varepsilon }(s) \gamma_{x}^{j}(s) \\
&&+\mathcal{S}(s) \delta \gamma_{j}^{\varepsilon }(s) +\mathcal{S}(s) \delta \gamma_{j}^{\varepsilon }(s) \gamma_{x}^{j}(s) .
\end{eqnarray*}

\begin{description}
\item From Lemma \ref{est}, we derive that
\begin{equation*}
\mathbb{E}\bigg[\bigg\vert g(x^{\varepsilon }(T))-g(\bar{x}(T)) -\mathfrak{p}(T)(x_{1}+x_{2})(T) -\frac{1}{2}\mathcal{P}(T)x_{1}(T)^{2}\bigg\vert^{2}\bigg] =o(\varepsilon^{2})
\end{equation*}
and
\begin{equation*}
\mathbb{E}\bigg[\bigg(\int_{0}^{T}\vert A_{4}(s) x_{1}(s) \vert \mathrm{d}s\bigg)^{2}\bigg] = o(\varepsilon^{2}).
\end{equation*}
Let
\begin{eqnarray*}
\bar{y}^{\varepsilon }(t) &=&y^{\varepsilon }(t) - \bigg[\mathfrak{p}(t)(x_{1}+x_{2})(t)-\frac{1}{2}\mathcal{P}(t) x_{1}(t)^{2}\bigg], \\
\bar{z}^{\varepsilon }(t) &=&z^{\varepsilon }(t) -\Big [\mathfrak{p}(s) \delta \sigma^{\varepsilon }(s) I_{E_{\varepsilon }}(s) +A_{5}(s)(x_{1}(s) +x_{2}(s)) \\
&&+\frac{1}{2}A_{6}(s) x_{1}(s) x_{1}(s)+A_{7}(s) x_{1}(s) I_{E_{\varepsilon }}(s) \Big ], \\
\bar{\kappa}_{j}^{\varepsilon }(t) &=&\kappa_{j}^{\varepsilon}(t) -\Big [(\mathfrak{p}(s) +\mathfrak{s}_{j}(s)) \delta \gamma_{j}^{\varepsilon }(s)I_{E_{\varepsilon }}(s) \\
&&+\frac{1}{2}(\mathcal{P}(s) +\mathcal{S}(s)) \delta \gamma_{j}^{\varepsilon }(s) \delta \gamma_{j}^{\varepsilon }(s) +A_{8}^{j}(s)(x_{1}(s) +x_{2}(s)) \\
&&+\frac{1}{2}A_{9}^{j}(s) x_{1}(s) x_{1}(s) +A_{10}^{j}(s) x_{1}(s) \Big].
\end{eqnarray*}
Then, a simple calculation yields
\begin{eqnarray*}
\bar{y}^{\varepsilon }(t)
&=&g(\bar{x}(T)) +\int_{t}^{T}\Big [f(s,x^{\varepsilon }(s),y^{\varepsilon }(s),z^{\varepsilon }(s),\kappa^{\varepsilon }(s),u^{\varepsilon }(s),\alpha(s-)) \\
&&+A_{1}(s) +A_{2}(s)(x_{1}(s)+x_{2}(s)) +\frac{1}{2}A_{3}(s) x_{1}(s)^{2}\Big ]\mathrm{d}s \\
&&-\int_{t}^{T}\bar{z}^{\varepsilon }(t) \mathrm{d}W(s) -\int_{t}^{T}\bar{\kappa}^{\varepsilon }(t) \mathrm{d}\tilde{\Phi}(s) +o(\varepsilon) .
\end{eqnarray*}
Now let
\begin{eqnarray*}
\hat{y}^{\varepsilon }(t) &=&\bar{y}^{\varepsilon }(t) -\bar{y}(t), \\
\hat{z}^{\varepsilon }(t) &=&\bar{z}^{\varepsilon }(t) -\bar{z}(t), \\
\hat{\kappa}^{\varepsilon }(t) &=&\bar{\kappa}^{\varepsilon}(t) -\bar{\kappa}(t) .
\end{eqnarray*}
It is easy to show that
\begin{eqnarray}
\hat{y}^{\varepsilon }(t) &=&o(\varepsilon)+\int_{t}^{T}\Big [f(s,x^{\varepsilon }(s),y^{\varepsilon}(s),z^{\varepsilon }(s),\kappa^{\varepsilon}(s),u^{\varepsilon }(s),\alpha(s-))  \notag \\
&&-f(s,\bar{x}(s),\bar{y}(s),\bar{z}(s),\bar{\kappa}(s),\bar{u}(s),\alpha(s-))  \notag \\
&&+A_{1}(s) +A_{2}(s)(x_{1}(s)+x_{2}(s)) +\frac{1}{2}A_{3}(s) x_{1}(s)^{2}\Big ]\mathrm{d}s  \notag \\
&&-\int_{t}^{T}\hat{z}^{\varepsilon }(s) \mathrm{d}W(s) -\int_{t}^{T}\hat{\kappa}^{\varepsilon }(s) \mathrm{d}\tilde{\Phi}(s) .  \label{verbsde1}
\end{eqnarray}
Next, we deal with
\begin{eqnarray*}
&&f(s,x^{\varepsilon }(s),y^{\varepsilon }(s),z^{\varepsilon }(s),\kappa^{\varepsilon }(s),u^{\varepsilon }(s),\alpha(s-)) \\
&&-f(s,\bar{x}(s),\bar{y}(s),\bar{z}(s),\bar{\kappa}(s),\bar{u}(s),\alpha(s-)) \\
&=&\Big [f\Big(s,\bar{x}(s),\bar{y}(s),\bar{z}(s) +\mathfrak{p}(s) \delta \sigma^{\varepsilon}(s), \\
&&\bar{\kappa}(s) +\Big [(\mathfrak{p}(s) +\mathfrak{s}_{j}(s)) \delta \gamma_{j}^{\varepsilon}(s) +\frac{1}{2}(\mathcal{P}(s) +\mathcal{S}(s)) \delta \gamma_{j}^{\varepsilon }(s) \delta \gamma_{j}^{\varepsilon }(s) \Big ]_{1\times D},u,\alpha(s-) \Big) \\
&&-f(s,\bar{x}(s),\bar{y}(s),\bar{z}(s),\bar{\kappa}(s),\bar{u}(s),\alpha(s-)) \Big ]I_{E_{\varepsilon }}(s) \\
&&+f(s,x^{\varepsilon }(s),y^{\varepsilon }(s),z^{\varepsilon }(s),\kappa^{\varepsilon }(s),u^{\varepsilon }(s),\alpha(s-)) \\
&&-f\Big(s,\bar{x}(s),\bar{y}(s),\bar{z}(s) +\mathfrak{p}(s) \delta \sigma^{\varepsilon }(s) I_{E_{\varepsilon }}(s),\bar{\kappa}(s) +\Big [(\mathfrak{p}(s) +\mathfrak{s}_{j}(s)) \delta \gamma_{j}^{\varepsilon }(s) \\
&&+\frac{1}{2}(\mathcal{P}(s) +\mathcal{S}(s)) \delta \gamma_{j}^{\varepsilon }(s) \delta \gamma_{j}^{\varepsilon }(s) \Big ]_{1\times D}I_{E_{\varepsilon}}(s),u,\alpha(s-) \Big) \\
&=&\Big [f\Big(s,\bar{x}(s),\bar{y}(s),\bar{z}(s) +\mathfrak{p}(s) \delta \sigma^{\varepsilon}(s), \\
&&\bar{\kappa}(s) +\Big [(\mathfrak{p}(s) +\mathfrak{s}_{j}(s)) \delta \gamma_{j}^{\varepsilon}(s) +\frac{1}{2}(\mathcal{P}(s) +\mathcal{S}(s)) \delta \gamma_{j}^{\varepsilon }(s)\delta \gamma_{j}^{\varepsilon }(s),u,\alpha(s-)\Big ]_{1\times D} \\
&&-f(s,\bar{x}(s),\bar{y}(s),\bar{z}(s),\bar{\kappa}(s),\bar{u}(s),\alpha(s-)) \Big ]I_{E_{\varepsilon }}(s) \\
&&+f(s,\bar{x}(s) +x_{1}(s) +x_{2}(s),\bar{y}(s) +\hat{y}^{\varepsilon }(s) +A_{11}(s), \\
&&\bar{z}(s) +\hat{z}^{\varepsilon }(s)+A_{12}(s),\bar{\kappa}(s) +\hat{\kappa}^{\varepsilon }(s) +A_{13}(s),\bar{u}(s),\alpha(s-)) \\
&&-f(s,\bar{x}(s),\bar{y}(s),\bar{z}(s),\bar{\kappa}(s),\bar{u}(s),\alpha(s-)),
\end{eqnarray*}
where
\begin{eqnarray*}
A_{11}(s) &=&\mathfrak{p}(t)(x_{1}+x_{2})(t) +\frac{1}{2}\mathcal{P}(t)x_{1}(t)^{2}, \\
A_{12}(s) &=&A_{5}(s)(x_{1}(s)+x_{2}(s)) +\frac{1}{2}A_{6}(s) x_{1}(s) x_{1}(s) +A_{7}(s) x_{1}(s)I_{E_{\varepsilon }}(s), \\
A_{13}(s) &=&\Big [A_{8}^{j}(s)(x_{1}(s) +x_{2}(s)) +\frac{1}{2}A_{9}^{j}(s)x_{1}(s) x_{1}(s) +A_{10}^{j}(s)x_{1}(s) I_{E_{\varepsilon }}(s) \Big ]_{1\times D}.
\end{eqnarray*}
We shall construct a linear BSDE with solution $(\hat{y}^{\varepsilon}(t),\hat{z}^{\varepsilon }(t),\hat{\kappa}^{\varepsilon }(t)) $ via selecting certain suitable $F(t)$, $G(t) $. The following conditions must be satisfied:
\end{description}

\begin{itemize}
\item $F(t) $ and $G(t) $ are determined by $(s,\bar{x}(s),\bar{y}(s),\bar{z}(s),\bar{\kappa}(s),\bar{u}(s),\alpha(s-))$;

\item We must have
\begin{eqnarray*}
&& f(s,\bar{x}(s) +x_{1}(s) +x_{2}(s),\bar{y}(s) +A_{11}(s),\bar{z}(s) +A_{12}(s),\bar{\kappa}(s)+A_{13}(s),\bar{u}(s),\alpha(s-)) \\
&& -f(s,\bar{x}(s),\bar{y}(s),\bar{z}(s),\bar{\kappa}(s),\bar{u}(s),\alpha(s-)) +A_{2}(s)(x_{1}(s) +x_{2}(s)) +\frac{1}{2}A_{3}(s) x_{1}(s)^{2} \\
&& = o(\varepsilon),
\end{eqnarray*}
where $o(\varepsilon)$ is independent on $x_{1}(s)$ and $x_{2}(s)$.
\end{itemize}

Applying Taylor's expansion to
\begin{eqnarray*}
&&f(s,\bar{x}(s) +x_{1}(s) +x_{2}(s),\bar{y}(s) +A_{11}(s),\bar{z}(s)+A_{12}(s),\bar{u}(s),\alpha(s-)) \\
&&-f(s,\bar{x}(s),\bar{y}(s),\bar{z}(s),\bar{\kappa}(s),\bar{u}(s),\alpha(s-)),
\end{eqnarray*}
we are able to find
\begin{eqnarray*}
F(s) &=& \bigg[b_{x}(s) +f_{y}(s)+f_{z}(s) \sigma_{x}(s) +\sum_{j=1}^{D}f_{\kappa}^{j}(s) \gamma_{x}^{j}(s)\bigg]\mathfrak{p}(s) \\
&&+ \big[\sigma_{x}(s) +f_{z}(s)\big]\mathfrak{q}(s) +\sum_{j=1}^{D}\Big(\gamma_{x}^{j}(s)\lambda_{j}(s) +f_{\kappa }^{j}(s)\Big)\mathfrak{s}_{j}(s) \\
&&+\sum_{j=1}^{D}f_{\kappa }^{j}(s) \delta \gamma_{x}^{j}(s) +f_{x}(t) \\
G(s) &=& \big[2b_{x}(s) +\sigma_{x}(s)\sigma_{x}(s) +\gamma_{x}(s) \gamma_{x}(s) \lambda_{j}(s) +f_{y}(s) +f_{z}(s) 2\sigma_{x}(s)\big]\mathcal{P}(s) \\
&&+\big[2\sigma_{x}(s) +f_{z}(s) \big]\mathcal{Q}(s) +b_{xx}(s) \mathfrak{p}(s) + \big[\mathfrak{q}(s) +f_{z}(s) \mathfrak{p}(s)\big]\sigma_{xx}(s) \\
&&+\sum_{j=1}^{D}\Big(\gamma_{xx}^{j}(s) \mathfrak{s}_{j}(s) +\mathcal{S}(s) \gamma_{x}^{j}(s) \gamma_{x}^{j}(s) +2\mathcal{S}(s) \gamma_{x}^{j}(s)\Big)\lambda_{j}(s) \\
&&+\sum_{j=1}^{D}f_{\kappa }^{j}(s)\Big(\mathfrak{p}(s)\gamma_{xx}^{j}(s) +\mathfrak{s}_{j}(s) \gamma_{xx}^{j}(s) +2\mathcal{P}(s) \gamma_{x}^{j}(s) +\mathcal{P}(s) \gamma_{x}^{j}(s) \gamma_{x}^{j}(s) \\
&&+\mathcal{S}(s) +2\mathcal{S}(s) \gamma_{x}^{j}(s) +\mathcal{S}(s) \gamma_{x}^{j}(s) \gamma_{x}^{j}(s)\Big) \\
&&+ \bigg(1,\mathfrak{p}(s),\mathfrak{p}(s) \sigma_{x}(s) +\mathfrak{q}(s),\sum_{j=1}^{D}(\mathfrak{p}(s) \gamma_{x}^{j}(s) +\mathfrak{s}_{j}(s) +\delta \gamma_{x}^{j}(s))\bigg)  \\
&&\cdot D^{2}f(t)\bigg(1,\mathfrak{p}(s),\mathfrak{p}(s) \sigma_{x}(s) +\mathfrak{q}(s),\sum_{j=1}^{D}(\mathfrak{p}(s) \gamma_{x}^{j}(s) +\mathfrak{s}_{j}(s) +\delta \gamma_{x}^{j}(s))\bigg)^{\top},
\end{eqnarray*}
where $f(t) =f(s,\bar{x}(s),\bar{y}(s),\bar{z}(s),\bar{\kappa}(s),\bar{u}(s),\alpha(s-)) $ and similarly for $f_{x}$, $f_{y}(\cdot)$, $f_{z}(\cdot)$, $f_{\kappa}(\cdot) $ and $D^{2}f(\cdot)$.
By Assumptions (A1) and(A2), we get that the adjoint equations (\ref{fbsdeadj1}) and (\ref{fbsdeadj2}) have unique solutions $(\mathfrak{p}(\cdot) \mathfrak{,q}(\cdot) \mathfrak{,s}(\cdot))$ and $(\mathcal{P}(\cdot),\mathcal{Q}(\cdot),\mathcal{S}(\cdot))$, respectively, and
\begin{equation*}
\mathbb{E}\bigg [\sup_{0\leq t\leq T}(\vert \mathfrak{p}(t) \vert^{2}+\vert \mathcal{P}(t) \vert^{2}) +\int_{0}^{T}\Big [\vert \mathfrak{q}(s)\vert^{2}+\vert \mathcal{Q}(s) \vert^{2}+\sum_{j=1}^{D}(\vert \mathfrak{s}_{j}(s)\vert^{2}+\vert \mathcal{S}_{j}(s) \vert^{2}) \Big ]\mathrm{d}s\bigg ]<\infty .
\end{equation*}
We introduce the following BSDE:
\begin{eqnarray}
\hat{y}(t) &=&\int_{t}^{T}\Big \{f_{y}(s) \hat{y}(s) +f_{z}(s) \hat{z}(s)+\sum_{j=1}^{D}f_{\kappa }^{j}(s) \hat{\kappa}_{j}(s) \lambda_{j}(t) \notag \\
&&+\Big [\mathfrak{p}(s)\delta b^{\varepsilon }(s)+\delta \sigma^{\varepsilon }(s) \mathfrak{q}(s) +\frac{1}{2}\mathcal{P}(s) \delta \sigma^{\varepsilon }(s) \delta \sigma^{\varepsilon }(s) \notag \\
&&+\sum_{j=1}^{D}\Big(\delta \gamma_{j}^{\varepsilon }(s)\mathfrak{s}_{j}(s) +\frac{1}{2}(\mathcal{P}(s) +\mathcal{S}(s)) \delta \gamma_{j}^{\varepsilon}(s) \delta \gamma_{j}^{\varepsilon }(s) \Big)\lambda_{j}(s) \notag \\
&&+f\Big(s,\bar{x}(s),\bar{y}(s),\bar{z}(s) +\mathfrak{p}(s) \delta \sigma^{\varepsilon }(s), \notag \notag \\
&&\bar{\kappa}(s) +\Big((\mathfrak{p}(s) +\mathfrak{s}_{j}(s)) \delta \gamma_{j}^{\varepsilon}(s) +\frac{1}{2}(\mathcal{P}(s) +\mathcal{S}(s)) \delta \gamma_{j}^{\varepsilon }(s)\delta \gamma_{j}^{\varepsilon }(s) \Big)_{1\times D},u,\alpha(s-) \Big) \notag \\
&&-f(s,\bar{x}(s),\bar{y}(s),\bar{z}(s),\bar{\kappa}(s),\bar{u}(s),\alpha(s-)) \Big ]I_{E_{\varepsilon }}(s) \Big \}\mathrm{d}s \notag \\
&&-\int_{t}^{T}\hat{z}(s) \mathrm{d}W(s) -\int_{t}^{T}\hat{\kappa}(s) \mathrm{d}\tilde{\Phi}(s) .
\label{verbsde2}
\end{eqnarray}

Next, we shall prove the following estimates:

\begin{lemma} \sl
Assume that \emph{(A1)-(A2)}
\begin{equation}
\mathbb{E}\bigg[\sup_{0\leq t\leq T}\vert \hat{y}^{\varepsilon }(t) \vert^{2}+\int_{0}^{T}\Big [\vert \hat{z}^{\varepsilon}(s) \vert^{2}+\sum_{j=1}^{D}\vert \hat{\kappa}_{j}^{\varepsilon }(s) \vert^{2}\lambda_{j}(s) \Big ]\mathrm{d}s\bigg] =O(\varepsilon^{2}),
\label{vest1}
\end{equation}
\begin{equation}
\mathbb{E}\bigg[\sup_{0\leq t\leq T}\vert \hat{y}(t)\vert^{2}+\int_{0}^{T}\Big [\vert \hat{z}(s)\vert^{2}+\sum_{j=1}^{D}\vert \hat{\kappa}_{j}(s)\vert^{2}\lambda_{j}(s) \Big ]\mathrm{d}s\bigg]=O(\varepsilon^{2}),  \label{vest2}
\end{equation}
\begin{eqnarray}
&&\mathbb{E}\bigg [\sup_{0\leq t\leq T}\vert \hat{y}^{\varepsilon}(t) -\hat{y}(t) \vert^{2}+\int_{0}^{T}\Big [\vert \hat{z}^{\varepsilon }(s) -\hat{z}(s)\vert^{2}  +\sum_{j=1}^{D}\vert \hat{\kappa}_{j}^{\varepsilon }(s) -\hat{\kappa}_{j}(s) \vert^{2}\lambda_{j}(s)\Big ]\mathrm{d}s\bigg] \notag \\
&& = o(\varepsilon^{2}) .  \label{vest3}
\end{eqnarray}
\end{lemma}

\paragraph{Proof.}

We first prove (\ref{vest1}). We reformulate BSDE (\ref{verbsde1}) as follows:
\begin{eqnarray*}
\Pi_{1}(s) & \!\!=\!\! & f(s,x^{\varepsilon }(s),y^{\varepsilon }(s),z^{\varepsilon }(s),\kappa^{\varepsilon }(s),u^{\varepsilon }(s),\alpha(s-)) \\
&&-f(s,(\bar{x}+x_{1}+x_{2})(s),(\bar{y}^{\varepsilon }+A_{11})(s), (\bar{z}^{\varepsilon }+A_{12})(s),(\bar{\kappa}^{\varepsilon }+A_{13})(s),\bar{u}(s),\alpha(s-)), \\
\Pi_{2}(s) & \!\!=\!\! & f(s,(\bar{x}+x_{1}+x_{2})(s),(\bar{y}^{\varepsilon }+A_{11})(s), (\bar{z}^{\varepsilon }+A_{12})(s),(\bar{\kappa}^{\varepsilon }+A_{13})(s),\bar{u}(s),\alpha(s-)) \\
&&-f((s,(\bar{x}+x_{1}+x_{2})(s),(\bar{y}+A_{11})(s), (\bar{z}+A_{12})(s),(\bar{\kappa}+A_{13})(s),\bar{u}(s),\alpha(s-)), \\
\Pi_{3}(s) & \!\!=\!\! & f(s,(\bar{x}+x_{1}+x_{2})(s),(\bar{y}+A_{11})(s), (\bar{z}+A_{12})(s),(\bar{\kappa}+A_{13})(s),\bar{u}(s),\alpha(s-)).
\end{eqnarray*}
Clearly,
\begin{eqnarray*}
&&f(s,x^{\varepsilon }(s),y^{\varepsilon }(s),z^{\varepsilon }(s),\kappa^{\varepsilon }(s),u^{\varepsilon }(s),\alpha(s-)) -f(s,\bar{x}(s),\bar{y}(s),\bar{z}(s),\bar{\kappa}(s),\bar{u}(s),\alpha(s-)) \\
&& = \Pi_{1}(s) +\Pi_{2}(s) +\Pi_{3}(s).
\end{eqnarray*}
First, we consider
\begin{equation*}
\Pi_{2}(s) =\check{f}_{y}(s) \hat{y}^{\varepsilon}(s) +\check{f}_{z}(s) \hat{z}^{\varepsilon }(s) +\check{f}_{\kappa }(s) \hat{\kappa}^{\varepsilon}(s),
\end{equation*}
where
\begin{eqnarray*}
\check{f}_{l}(s) &=&\int_{0}^{1}f_{l}(s,(\bar{x}+x_{1}+x_{2})(s),\theta((\bar{y}^{\varepsilon }+A_{11})(s)) +(1-\theta)(\bar{y}+A_{11})(s), \\
&&\theta((\bar{z}^{\varepsilon }+A_{12})(s)) +(1-\theta)(\bar{z}+A_{12})(s), \\
&&\theta((\bar{\kappa}^{\varepsilon }+A_{13})(s)) +(1-\theta)(\bar{\kappa}+A_{13})(s),\bar{u}(s),\alpha(s-))\mathrm{d}\theta,
\end{eqnarray*}
for $l=y$, $z$, $\kappa$, respectively.
Next, we deal with the following term
\begin{eqnarray*}
\Pi_{3}(s) &=&f_{x}(s)(x_{1}+x_{2})(s) +f_{y}(s) A_{11}(s) +f_{z}(s) A_{12}(s) +f_{\kappa }(s) A_{13}(s) \\
&&+\frac{1}{2}\big((x_{1}+x_{2})(s),A_{11}(s),A_{12}(s),A_{13}(s)\big) \\
&& \cdot\check{D}^{2}(s)\big((x_{1}+x_{2})(s),A_{11}(s),A_{12}(s),A_{13}(s)\big)^{\top},
\end{eqnarray*}
where
\begin{eqnarray*}
\check{D}^{2}(s) &=& 2\int_{0}^{1}\int_{0}^{1}\theta D^{2}f\big(s,\bar{x}(s) +\theta \mu(x_{1}+x_{2})(s),\bar{y}(s) +\theta \mu A_{11}(s), \\
&& \quad \bar{z}(s) +\theta \mu A_{12}(s),\bar{\kappa}(s) +\theta \mu A_{13}(s),\bar{u}(s),\alpha(s-)\big)\mathrm{d}\theta \mathrm{d}\mu.
\end{eqnarray*}
It therefore follows that
\begin{eqnarray}
\hat{y}^{\varepsilon }(t) &\!\!=\!\!&o(\varepsilon)+\int_{t}^{T}\bigg [A_{1}(s) I_{E_{\varepsilon }}(s)+\Pi_{1}(s) +\check{f}_{y}(s) \hat{y}^{\varepsilon}(s) +\check{f}_{z}(s) \hat{z}^{\varepsilon }(s) +\check{f}_{\kappa }(s) \hat{\kappa}^{\varepsilon}(s)  \notag \\
&&\frac{1}{2}\big((x_{1}+x_{2})(s),A_{11}(s),A_{12}(s),A_{13}(s)\big)\check{D}^{2}(s)\big((x_{1}+x_{2})(s),A_{11}(s),A_{12}(s),A_{13}(s)\big)^{\top}  \notag \\
&&+\Big [f_{z}(s) A_{7}(s) x_{1}(s)+\sum_{j=1}^{D}A_{10}^{j}(s) x_{1}(s) \Big]I_{E_{\varepsilon}}(s) -\frac{1}{2}\Pi_{4}(s) D^{2}f(s) \Pi_{4}(s)^{\top }\bigg]\mathrm{d}s \notag \\
&&-\int_{t}^{T}\hat{z}^{\varepsilon}(s) \mathrm{d}W(s) -\int_{t}^{T}\hat{\kappa}^{\varepsilon}(s) \mathrm{d}\tilde{\Phi}(s), \label{verbsde3}
\end{eqnarray}
where
\begin{equation*}
\Pi_{4}(s) = \Big(1,\mathfrak{p}(s),\mathfrak{p}(s) \sigma_{x}(s) +\mathfrak{q}(s),\sum_{j=1}^{D}(\mathfrak{p}(s) \gamma_{x}^{j}(s) +\mathfrak{s}_{j}(s) +\delta \gamma_{x}^{j}(s))\Big).
\end{equation*}
It follows from Lemma \ref{baest} and Lemma \ref{est} that we obtain (\ref{vest1}). By the same argument, we can easily get (\ref{vest2}). Next we handle (\ref{vest3}). To this end, we set
\begin{eqnarray*}
\tilde{x}^{\varepsilon }(t) &=&x^{\varepsilon }(t) -\bar{x}(t) -x_{1}(t) -x_{2}(t), \\
\tilde{y}^{\varepsilon }(t) &=&\hat{y}^{\varepsilon }(t) -\hat{y}(t), \\
\tilde{z}^{\varepsilon }(t) &=&\hat{z}^{\varepsilon }(t) -\hat{z}(t), \\
\tilde{\kappa}^{\varepsilon }(t) &=&\hat{\kappa}^{\varepsilon}(t) -\hat{\kappa}(t) .
\end{eqnarray*}
Subtracting (\ref{verbsde2}) from BSDE (\ref{verbsde3}) yields
\begin{eqnarray}
&& \tilde{y}^{\varepsilon }(t) = o(\varepsilon) \notag \\
&& +\int_{t}^{T}\Big [\tilde{f}_{y}(s) \tilde{y}^{\varepsilon}(s) +\tilde{f}_{z}(s) \tilde{z}^{\varepsilon}(s) +\check{f}_{\kappa }(s) \tilde{\kappa}^{\varepsilon }(s) +(\tilde{f}_{y}-f_{y})(s) \hat{y}(s) \notag \\
&& +(\tilde{f}_{z}-f_{z})(s) \hat{z}(s)+(\tilde{f}_{\kappa }-f_{\kappa })(s) \hat{\kappa}(s) +\Pi_{1}(s)  \notag \\
&& -\Big [f\Big(s,\bar{x}(s),\bar{y}(s),\bar{z}(s) +\mathfrak{p}(s) \delta\sigma^{\varepsilon }(s), \notag \\
&& \bar{\kappa}(s) +\Big((\mathfrak{p}(s) +\mathfrak{s}_{j}(s)) \delta \gamma_{j}^{\varepsilon}(s) +\frac{1}{2}(\mathcal{P}(s) +\mathcal{S}(s)) \delta \gamma_{j}^{\varepsilon }(s) \delta \gamma_{j}^{\varepsilon }(s) \Big)_{1\times D},u,\alpha(s-) \Big) \notag \\
&& -f(s) \Big ]I_{E_{\varepsilon }}(s) +\Big [f_{z}(s) A_{7}(s) x_{1}(s)+\sum_{j=1}^{D}A_{10}^{j}(s) x_{1}(s) \Big ]I_{E_{\varepsilon }}(s)  \notag \\
&& +\frac{1}{2}\big((x_{1}+x_{2})(s),A_{11}(s),A_{12}(s),A_{13}(s)\big)\check{D}^{2}f(s)\big((x_{1}+x_{2})(s),A_{11}(s),A_{12}(s),A_{13}(s)\big)^{\top}  \notag \\
&& +\Big [f_{z}(s) A_{7}(s)+\sum_{j=1}^{D}A_{10}^{j}(s) \Big ]x_{1}(s)I_{E_{\varepsilon }}(s) -\frac{1}{2}\Pi_{4}(s) D^{2}f(s) \Pi_{4}(s)^{\top }x_{1}(s) x_{1}(s) \Big ]\mathrm{d}s \notag \\
&& -\int_{t}^{T}\hat{z}^{\varepsilon }(s) \mathrm{d}W(s) -\int_{t}^{T}\hat{\kappa}^{\varepsilon }(s) \mathrm{d}\tilde{\Phi}(s) .  \label{verbsde4}
\end{eqnarray}
Applying Lemma \ref{baest}, we only need to prove that
\begin{equation}
\mathbb{E}\bigg[\bigg(\int_{0}^{T}\vert(\tilde{f}_{y}-f_{y})(s) \hat{y}(s) +(\tilde{f}_{z}-f_{z})(s) \hat{z}(s) +(\tilde{f}_{\kappa }-f_{\kappa })(s) \hat{\kappa}(s)\vert\bigg)^{2}\mathrm{d}s\bigg] =o(\varepsilon^{2}),
\label{k1}
\end{equation}
\begin{eqnarray}
&& \mathbb{E}\bigg [\bigg(\int_{0}^{T}\Big\vert\Pi_{1}(s) -\Big[f\Big(s,\bar{x}(s),\bar{y}(s),\bar{z}(s) +\mathfrak{p}(s) \delta \sigma^{\varepsilon }(s),  \notag \\
&& \bar{\kappa}(s) +\Big((\mathfrak{p}(s) +\mathfrak{s}_{j}(s)) \delta \gamma_{j}^{\varepsilon}(s) +\frac{1}{2}(\mathcal{P}(s) +\mathcal{S}(s)) \delta \gamma_{j}^{\varepsilon }(s)\delta \gamma_{j}^{\varepsilon }(s) \Big)_{1\times D},u,\alpha(s-) \Big) \notag \\
&& -f(s) \Big ]I_{E_{\varepsilon }}(s) \Big\vert\mathrm{d}s\bigg)^{2}\bigg ]= o(\varepsilon^{2}),  \label{k2}
\end{eqnarray}
and
\begin{equation}
\mathbb{E}\bigg[\bigg(\int_{0}^{T}\Big\vert \Pi_{4}(s)\big[\tilde{D}^{2}f(s) -D^{2}f(s)\big]\Pi_{4}(s)^{\top}x_{1}(s) x_{1}(s)\Big\vert\bigg)^{2}\mathrm{d}s\bigg] = o(\varepsilon^{2}).  \label{k3}
\end{equation}
Obviously,
\begin{eqnarray}
&&\vert(\tilde{f}_{y}-f_{y})(s) + \tilde{f}_{z}-f_{z})(s) +(\tilde{f}_{\kappa}-f_{\kappa })(s) \vert  \label{k11} \\
& \leq &C(\vert(x_{1}+x_{2})(s)\vert +\vert A_{11}(s) \vert +\vert
A_{12}(s) \vert +\vert A_{13}(s)\vert +\vert \hat{y}^{\varepsilon }(s) \vert+\vert \hat{z}^{\varepsilon }(s) \vert +\vert\hat{\kappa}^{\varepsilon }(s) \vert) .  \notag
\end{eqnarray}
One needs to check that
\begin{eqnarray}
\mathbb{E}\bigg[\bigg(\int_{0}^{T}\vert \mathfrak{q}(s)x_{1}(s) \hat{z}(s) \vert \mathrm{d}s\bigg)^{2}\bigg] &\leq &\mathbb{E}\bigg[\sup_{0\leq t\leq T}\vert x_{1}(t) \vert^{2}\bigg(\int_{0}^{T}\vert \mathfrak{q}(s) \vert^{2}\mathrm{d}s\bigg)\bigg(\int_{0}^{T}\vert \hat{z}(s) \vert^{2}\mathrm{d}s\bigg)\bigg]  \notag \\
&=&o(\varepsilon^{2}),  \label{k22}
\end{eqnarray}
and
\begin{eqnarray}
&&\mathbb{E}\bigg[\bigg(\int_{0}^{T}\sum_{j=1}^{D}\vert \mathfrak{s}_{j}(s) x_{1}(s) \hat{\kappa}_{j}(s)\lambda_{j}(s) \vert \mathrm{d}s\bigg)^{2}\bigg] \notag \\
&\leq &\mathbb{E}\bigg[\sup_{0\leq t\leq T}\vert x_{1}(t)\vert^{2}\sum_{j=1}^{D}\Big(\int_{0}^{T}\vert \mathfrak{s}_{j}(s) \vert^{2}\mathrm{d}s\Big)\Big(\int_{0}^{T}\vert \hat{\kappa}_{j}(s) \lambda_{j}(s) \vert^{2}\mathrm{d}s\Big)\bigg]  \notag \\
&=& o(\varepsilon^{2}) .  \label{k33}
\end{eqnarray}
From (\ref{k11})-(\ref{k33}), we prove (\ref{k1}). Due to $D^{2}f$ is
bounded the definition of $\check{D}^{2}f(s)$, we have
\begin{equation*}
\lim_{\varepsilon \rightarrow 0}\mathbb{E}\bigg[\bigg(\int_{0}^{T}\vert \Pi_{4}(s)\big[\tilde{D}^{2}f(s) -D^{2}f(s)\big]\Pi_{4}(s)^{\top}\vert\bigg)^{2}\mathrm{d}s\bigg] = 0,
\end{equation*}
which yields
\begin{eqnarray}
&&\mathbb{E}\bigg[\bigg(\int_{0}^{T}\vert \Pi_{4}(s)[ \tilde{D}^{2}f(s) -D^{2}f(s) ] \Pi_{4}(s)^{\top }x_{1}(s) x_{1}(s)\vert\bigg)^{2}\mathrm{d}s\bigg]  \notag \\
&\leq &\mathbb{E}[ \sup_{0\leq t\leq T}\vert x_{1}(s)\vert^{4}\int_{0}^{T}\vert \Pi_{4}(s) [ \tilde{D}^{2}f(s) -D^{2}f(s) ] \Pi_{4}(s)^{\top }\vert^{2}\mathrm{d}s]  \notag \\
&\leq &\varepsilon^{2}\mathbb{E}[ \int_{0}^{T}\vert \Pi_{4}(s) [ \tilde{D}^{2}f(s) -D^{2}f(s) ] \Pi_{4}(s)^{\top }\vert^{2}\mathrm{d}s]  \notag \\
&=&o(\varepsilon^{2}).  \label{k44}
\end{eqnarray}
It is fairly easy to get that
\begin{eqnarray*}
&& \Pi_{1}(s) -\Big [f\Big(s,\bar{x}(s),\bar{y}(s),\bar{z}(s) +\mathfrak{p}(s) \delta\sigma^{\varepsilon }(s), \\
&& \bar{\kappa}(s) +\sum_{j=1}^{D}\big((\mathfrak{p}(s) +\mathfrak{s}_{j}(s)) \delta \gamma_{j}^{\varepsilon }(s) +\frac{1}{2}(\mathcal{P}(s) +\mathcal{S}(s)) \delta \gamma_{j}^{\varepsilon }(s) \delta \gamma_{j}^{\varepsilon }(s) \big),u,\alpha(s-) \big) \\
&&  -f(s) \Big]I_{E_{\varepsilon }}(s) \\
&=& f(s,x^{\varepsilon }(s),y^{\varepsilon }(s),z^{\varepsilon }(s),\kappa^{\varepsilon }(s),u^{\varepsilon }(s),\alpha(s-))-f\Big(s,(\bar{x}+x_{1}+x_{2})(s),(\bar{y}^{\varepsilon }+A_{11})(s), \\
&& (\bar{z}^{\varepsilon }+A_{12}+\mathfrak{p}(s) \delta\sigma^{\varepsilon }(s) I_{E_{\varepsilon }(s)})(s),\big(\bar{\kappa}^{\varepsilon }+A_{13}+\big[(\mathfrak{p}(s) +\mathfrak{s}_{j}(s))
\delta \gamma_{j}^{\varepsilon }(s) \\
&& +\frac{1}{2}(\mathcal{P}(s) +\mathcal{S}(s)) \delta \gamma_{j}^{\varepsilon }(s) \delta \gamma_{j}^{\varepsilon }(s)\big]_{1\times D}I_{E_{\varepsilon}(s) }\big)(s),u^{\varepsilon }(s),\alpha(s-)\Big) \\
&&+\Big \{f\Big(s,(\bar{x}+x_{1}+x_{2})(s),(\bar{y}^{\varepsilon }+A_{11})(s), \\
&&(\bar{z}^{\varepsilon }+A_{12}+\mathfrak{p}(s) \delta\sigma^{\varepsilon }(s))(s),\big(\bar{\kappa}^{\varepsilon }+A_{13}(s) +\big[(\mathfrak{p}(s) +\mathfrak{s}_{j}(s)) \delta \gamma_{j}^{\varepsilon }(s) \\
&&+\frac{1}{2}(\mathcal{P}(s) +\mathcal{S}(s)) \delta \gamma_{j}^{\varepsilon }(s) \delta \gamma_{j}^{\varepsilon }(s)\big]_{1\times D}\big)(s),u^{\varepsilon }(s),\alpha(s-)\Big) \\
&&-f\big(s,(\bar{x}+x_{1}+x_{2})(s),(\bar{y}+A_{11})(s),(\bar{z}+A_{12})(s),(\bar{\kappa}+A_{13})(s),\bar{u}(s),\alpha(s-)\big) \\
&&-f\Big(s,\bar{x}(s),\bar{y}(s),\bar{z}(s) +\mathfrak{p}(s) \delta \sigma^{\varepsilon}(s), \\
&&\bar{\kappa}(s) +\Big [(\mathfrak{p}(s) +\mathfrak{s}_{j}(s)) \delta \gamma_{j}^{\varepsilon}(s) +\frac{1}{2}(\mathcal{P}(s) +\mathcal{S}(s)) \delta \gamma_{j}^{\varepsilon }(s)\delta \gamma_{j}^{\varepsilon }(s) \Big ]_{1\times D},u^{\varepsilon }(s),\alpha(s-) \Big)\\
&&+f\big(s,\bar{x}(s),\bar{y}(s),\bar{z}(s),\bar{\kappa}(s),\bar{u}(s),\alpha(s-)\big)\Big\}I_{E_{\varepsilon }(s) } \\
&\leq &C\Big \{\vert \tilde{x}^{\varepsilon }(s)\vert +\Big(\vert x_{1}(s) +x_{2}(s)\vert +\vert \hat{y}^{\varepsilon }(s) \vert+\vert \hat{z}^{\varepsilon }(s) \vert +\vert\hat{\kappa}^{\varepsilon }(s) \vert \\
&&+\vert A_{11}(s) \vert +\vert A_{12}(s) \vert +\vert A_{13}(s) \vert \Big)I_{E_{\varepsilon }(s) }\Big \}.
\end{eqnarray*}
One can check that
\begin{eqnarray*}
\mathbb{E}\bigg[\bigg(\int_{0}^{T}\vert \mathfrak{q}(s)x_{1}(s) \vert I_{E_{\varepsilon }}(s) \mathrm{d}s\bigg)^{2}\bigg]
&\leq &\varepsilon \mathbb{E}\bigg[\sup_{0\leq t\leq T}\vert x_{1}(t) \vert^{2}\cdot \int_{E_{\varepsilon}}\vert \mathfrak{q}(s) \vert^{2}\mathrm{d}s\bigg] \\
&\leq &\varepsilon^{2}\mathbb{E}\bigg[\int_{E_{\varepsilon }}\vert\mathfrak{q}(s) \vert^{2}\mathrm{d}s\bigg] \\
&=&o(\varepsilon^{2})
\end{eqnarray*}
and
\begin{equation*}
\mathbb{E}\bigg[\bigg(\int_{0}^{T}\vert \mathfrak{s}_{j}(s) x_{1}(s) \vert \lambda_{j}(s)I_{E_{\varepsilon }}(s) \mathrm{d}s\bigg)^{2}\bigg] =o(\varepsilon^{2}).
\end{equation*}
Therefore, (\ref{k3}) holds. Finally, from Lemma \ref{baest} we complete the proof.
\hfill $\Box $

Now we are able to show the variational equations for BSDE (\ref{bsde1}):
\begin{eqnarray}
\hat{y}^{\varepsilon }(t) &=&\bar{y}(t) +\mathfrak{p}(t)(x_{1}+x_{2})(t) +\frac{1}{2}\mathcal{P}(t) x_{1}(t) x_{1}(t) +\hat{y}(t) +o(\varepsilon),  \label{veq1} \\
\hat{z}^{\varepsilon }(t) &=&\bar{z}(t) +\mathfrak{p}(s) \delta \sigma^{\varepsilon }(s)I_{E_{\varepsilon }}(s) +A_{5}(s)(x_{1}(s) +x_{2}(s))  \notag \\
&&+\frac{1}{2}A_{6}(s) x_{1}(s) x_{1}(s)+A_{7}(s) x_{1}(s) I_{E_{\varepsilon }}(s) +\hat{z}(t) +o(\varepsilon),  \label{veq2} \\
\hat{\kappa}^{\varepsilon }(t) &=&\bar{\kappa}(t) +\Big [(\mathfrak{p}(s) +\mathfrak{s}_{j}(s)) \delta \gamma_{j}^{\varepsilon }(s) I_{E_{\varepsilon}}(s)  \notag \\
&&+\frac{1}{2}(\mathcal{P}(s) +\mathcal{S}(s)) \delta \gamma_{j}^{\varepsilon }(s) \delta \gamma_{j}^{\varepsilon }(s) +A_{8}^{j}(s)(x_{1}(s) +x_{2}(s))  \notag \\
&&+\frac{1}{2}A_{9}^{j}(s) x_{1}(s) x_{1}(s) +A_{10}^{j}(s) x_{1}(s) \Big ]_{1\times D}+\hat{\kappa}(t) +o(\varepsilon) .  \label{veq3}
\end{eqnarray}
Similar to $x_{1}(\cdot) $ and $x_{2}(\cdot)$, we have
\begin{equation*}
\left\{\begin{array}{rcl}
y_{1}(t) & = & \mathfrak{p}(t) x_{1}(t), \\
z_{1}(t) & = & \mathfrak{p}(s) \delta \sigma^{\varepsilon }(s) I_{E_{\varepsilon }}(s)+A_{5}(s) x_{1}(s), \\
\kappa_{1}(t) & = & \Big[(\mathfrak{p}(s) +\mathfrak{s}_{j}(s)) \delta \gamma_{j}^{\varepsilon}(s) I_{E_{\varepsilon }}(s) \\
&  & +\displaystyle\frac{1}{2}(\mathcal{P}(s) +\mathcal{S}(s)) \delta \gamma_{j}^{\varepsilon }(s) \delta\gamma_{j}^{\varepsilon }(s) I_{E_{\varepsilon }}(s) +A_{8}^{j}(s) x_{1}(s)\Big]_{1\times D}
\end{array}\right.
\end{equation*}
and
\begin{equation*}
\left\{\begin{array}{rcl}
y_{2}(t) & = & \mathfrak{p}(t)(x_{1}+x_{2})(t) +\displaystyle\frac{1}{2}\mathcal{P}(t)x_{1}(t) x_{1}(t) +\hat{y}(t), \\
z_{2}(t) & = & A_{5}(s)(x_{1}(s)+x_{2}(s)) +\displaystyle\frac{1}{2}A_{6}(s) x_{1}(s) x_{1}(s) +A_{7}(s) x_{1}(s) I_{E_{\varepsilon }}(s) +\hat{z}(t), \\
\kappa_{2}(t) & = & \Big[A_{8}^{j}(s) x_{2}(s) +\displaystyle\frac{1}{2}A_{9}^{j}(s) x_{1}(s)x_{1}(s) +A_{10}^{j}(s) x_{1}(s) I_{E_{\varepsilon}}(s)\Big]_{1\times D}+\hat{\kappa}(t).
\end{array}\right.
\end{equation*}
Now let us come back to discuss about the maximum principle for optimal control of FBSDEs (\ref{sde1})-(\ref{bsde1}). From the definition of performance functional (\ref{costfbsde}), it follows
\begin{eqnarray*}
J(u^{\varepsilon }(\cdot)) -J(\bar{u}(\cdot)) &=&y^{\varepsilon }(0) -\bar{y}(0) \\
&=&\hat{y}(0) +o(\varepsilon) \\
&\geq &0.
\end{eqnarray*}
Define the adjoint equation for (\ref{bsde1}) ($\lambda_{j}\neq 0$):
\begin{equation}
\left\{\begin{array}{rcl}
\mathrm{d}\chi(t) & = & f_{y}(t) \chi(t) \mathrm{d}t+f_{z}(t) \chi(t) \mathrm{d}W(t) +\displaystyle\sum_{j=1}^{D}\frac{1}{\lambda_{j}(t) }f_{\kappa }^{j}(t) \chi(t) \mathrm{d}\tilde{\Phi}_{j}(t), \\
\chi(0) & = & 1.
\end{array}\right.  \label{adjsde}
\end{equation}
Clearly, SDE (\ref{adjsde}) admits a unique strong solution. By virtue of It\^{o}'s formula to $\hat{y}(t) \chi(t) $, we can obtain
\begin{eqnarray*}
\hat{y}(0) &=&\mathbb{E}\bigg [\int_{0}^{T}\chi(s)\Big \{\mathfrak{p}(s) \delta b^{\varepsilon }(s)+\delta \sigma^{\varepsilon }(s) \mathfrak{q}(s) +\frac{1}{2}\mathcal{P}(s) \delta \sigma^{\varepsilon }(s) \delta \sigma^{\varepsilon }(s) \\
&&+\sum_{j=1}^{D}\Big [\delta \gamma_{j}^{\varepsilon }(s)
\mathfrak{s}_{j}(s) +\frac{1}{2}(\mathcal{P}(s) +\mathcal{S}(s)) \delta \gamma_{j}^{\varepsilon}(s) \delta \gamma_{j}^{\varepsilon }(s) \Big]\lambda_{j}(s) \\
&&+f\Big(s,\bar{x}(s),\bar{y}(s),\bar{z}(s) +\mathfrak{p}(s) \delta \sigma^{\varepsilon }(s), \\
&&\bar{\kappa}(s) +\Big [(\mathfrak{p}(s) +\mathfrak{s}_{j}(s)) \delta \gamma_{j}^{\varepsilon}(s)+\frac{1}{2}(\mathcal{P}(s) +\mathcal{S}(s)) \delta \gamma_{j}^{\varepsilon }(s) \delta \gamma_{j}^{\varepsilon }(s) \Big ]_{1\times D},u,\alpha(s-) \Big) \\
&&-f(s,\bar{x}(s),\bar{y}(s),\bar{z}(s),\bar{\kappa}(s),\bar{u}(s),\alpha(s-)) \Big \}I_{E_{\varepsilon }}(s) \mathrm{d}s\bigg].
\end{eqnarray*}
In order to state our main result, we define a function as follows:
\begin{eqnarray}
&&\mathbb{H}(t,x,y,z,\kappa,u,e_{i},\mathfrak{p,q,s},\mathcal{P},\mathcal{S})  \notag \\
&=&\mathfrak{p}b(t,x,u,e_{i}) +\sigma(t,x,u,e_{i})\mathfrak{q}+\frac{1}{2}\mathcal{P}(\sigma(t,x,u,e_{i})-\sigma(t,\bar{x},\bar{u},e_{i}))^{2}  \notag \\
&&+\sum_{j=1}^{D}\Big [\gamma_{j}(t,x,u,e_{i}) \mathfrak{s}+\frac{1}{2}(\mathcal{P}+\mathcal{S})(\gamma_{j}(t,x,u,e_{i}) -\gamma_{j}(t,\bar{x},\bar{u},e_{i}))^{2}\Big ]\lambda_{j}(t)  \notag \\
&&+f\Big(s,x,y,z+\mathfrak{p}(\sigma(t,x,u,e_{i})-\sigma(t,\bar{x},\bar{u},e_{i})),  \notag \\
&&\kappa +\Big [(\mathfrak{p}(s) +\mathfrak{s}_{j}(s))(\gamma(t,x,u,e_{i}) -\gamma(t,\bar{x},\bar{u},e_{i}))  \notag \\
&&+\frac{1}{2}(\mathcal{P}(s) +\mathcal{S}(s))(\gamma_{j}(t,x,u,e_{i}) -\gamma_{j}(t,\bar{x},\bar{u},e_{i}))^{2}\Big ]_{1\times D},u,e_{i}\Big).
\label{HM}
\end{eqnarray}
We now assert the following main result.

\begin{theorem}\label{the2} \sl
Assume that \emph{(A1)-(A2)} hold. Let $\bar{u}(\cdot)$ be an optimal control and $(\bar{x}(\cdot),\bar{y}(\cdot)$, $\bar{z}(\cdot),\bar{\kappa}(\cdot))$ be the associated solution. Then
\begin{eqnarray*}
&&\mathbb{H}(t,\bar{x}(t),\bar{y}(t),\bar{z}(t),\bar{\kappa}(t),u,\alpha(t),\mathfrak{p(t),q(t),s}(t),\mathcal{P}(t),\mathcal{S}(t)) \\
&\geq &\mathbb{H}(t,\bar{x}(t),\bar{y}(t),\bar{z}(t),\bar{\kappa}(t),\bar{u},\alpha(t),\mathfrak{p(t),q(t),s}(t),\mathcal{P}(t),\mathcal{S}(t)),  \quad \forall u\in U,\text{ a.e., a.s.,}
\end{eqnarray*}
where $\mathbb{H}$ is defined in (\ref{HM}).
\end{theorem}

We also provide a concrete example.

\begin{example}
Consider the following FBSDEs with $d=1$, $L=1$
\begin{equation*}
\left\{\begin{array}{rcl}
\mathrm{d}x(t) & = & u(t) \mathrm{d}W(t)+u(t)\displaystyle\sum_{j=1}^{D}\nu_{j}(t) \mathrm{d}\tilde{\Phi}_{j}(t), \\
-\mathrm{d}y(t) & = & f(z(t),\kappa(s)) \mathrm{d}t-z(t) \mathrm{d}W(t)-\kappa(t) \mathrm{d}\tilde{\Phi}(t), \\
x(0) & = & 0,\text{ }y(T) =x(T) +\alpha(T)^{\top }A\alpha(T),
\end{array}\right.
\end{equation*}
where $A\in \mathbb{R}^{D\times D}$ and $\nu_{j}(t) $ is a $\mathbb{R}$-valued deterministic process, $1\leq j\leq D$. Clearly, the
solutions to equations (\ref{fbsdeadj1}) and (\ref{fbsdeadj2}) are $(\mathfrak{p}(t) \mathfrak{,q(t),s}(t)) =(1,0,0)$, $(\mathcal{P}(t),\mathcal{Q}(t),\mathcal{S}(t)) =(0,0,0)$, respectively. Therefore, from Theorem \ref{the2}, the necessary condition for optimal control is
\begin{equation*}
f\bigg(\bar{z}(t)+u-\bar{u}(t),\bar{\kappa}(t) +(u-\bar{u}(t)) \sum_{j=1}^{D}\nu_{j}(t)\bigg) -f(\bar{z}(t),\bar{\kappa}(t)) \geq 0.
\end{equation*}
For instance, let $f(z,\kappa) =z+\sum_{j=1}^{D}\vert\kappa_{j}\vert^{2}\lambda_{j}$. Then, we immediately derive that $\bar{u}\equiv 0$ is an optimal control. The corresponding trajectories are $(0,0,0,0)$.
\end{example}

\section{Concluding remarks}

\label{Remarks}

In this paper, we study a general stochastic maximum principle for optimal
control for systems governed by a continuous-time Markov regime-switching
stochastic recursive utility. The control domain is assumed not to be
convex, and the diffusion terms allow to contain control variables.
On the one hand, we derive a result for forward stochastic optimal control
problems. Afterwards, based on previous estimates and inspirations, we
introduce two groups of new first and second-order adjoint equations.
Moreover, the corresponding variational equations for forward-backward
stochastic differential equations are obtained. Finally, we present
our results using an illustrated example. It is necessary to point out that,
due to the Markov regime-switching, the generator in the maximum principle
involves solutions of the novel second-order adjoint equations.
However, a number of other issues deserve further investigations as follows:
\begin{itemize}
\item Apart from Pontryagin's maximum principle, the Bellman's dynamic
programming principle (DPP for short) is another important tool in
tackling stochastic optimal control problems. Nowadays, many researchers
engage in this field and achieve fundamental result. As for the
DPP, systematic investigations of the classical stochastic optimal control
problems are discussed deeply in the famous book by Fleming and Soner \cite{FS}. Now
return to the expression of HJB equation (\ref{HJB}). It vividly indicates
that it is a coupled system of nonlinear differential equations because of
the presence of regime-switching, in contrast to a single nonlinear
differential equation in the traditional literature such as the user's guide
by Crandall et al \cite{CIL}. Out of question, it is
conceivable that the fruitful analysis will be much more involved than that
in \cite{CIL}. Nevertheless, we are not able to provide the
explicit form of the value function and an optimal control via solving (\ref{HJB}). The main obstacle appears since it is a coupled system of nonlinear
second-order partial differential equations, rendering extreme difficulty in
constructing a closed-form solution of (\ref{HJB}). In future, we shall
employ the viscosity solution to study this problem(see \cite{SZ, ZC} in this direction).

\item An interesting topic is to consider the fully coupled FBSDEs. Indeed,
some financial optimization problems for large investors and some asset
pricing problems with forward-backward differential utility (see \cite{CM,CC, RH} directly lead to fully coupled FBSDEs. Using the ideas from Peng and Wu \cite{PW}, under certain $G$-monotonicity conditions for
the coefficients, we may derive the existence and uniqueness of fully
coupled FBSDEs with Markov regime-switching. And then, a series of questions
turn out, for instance, the associated partial differential equations and
control problems. Another problem of great interests is to consider the case
when the random environment or the Markov regime-switching $\alpha$ is unobservable.
\end{itemize}

\appendix\label{app}

\section{The Proof of Lemma \protect\ref{baest}}

\paragraph{Proof of Lemma \protect\ref{baest}.}

We first assume that $F$ are bounded. Let $\epsilon >0$ and define
\begin{equation*}
\langle y\rangle_{\epsilon }\triangleq(\vert y\vert^{2}+\epsilon^{2})^{\frac{1}{2}},  \quad y\in \mathbb{R}.
\end{equation*}
Clearly for any $\epsilon >0$, the function $y \rightarrow \langle y\rangle_{\epsilon }$ is smooth
and $\langle y\rangle_{\epsilon} \rightarrow \vert y\vert $ as $\epsilon \rightarrow 0$.
The motivation of introducing such a function is to avoid some
difficulties that might be encountered in differentiating functions like $\vert y\vert^{2k}$ for non-integer $k$. Applying It\^{o}'s
formula to $\langle y\rangle_{\epsilon }^{2k}$ on $[t,T]$, we have
\begin{eqnarray}
&&\mathbb{E}[ \langle y(t) \rangle_{\epsilon}^{2k}] +k(2k-1) \mathbb{E}\int_{t}^{T}\langle y(s) \rangle_{\epsilon }^{2k-2}\big[\vert z(s) \vert^{2}+\Vert \kappa(s)\Vert\big] \mathrm{d}s  \notag \\
&\leq &\mathbb{E}[ \langle \xi \rangle_{\epsilon }^{2k}] +2k\mathbb{E}\bigg[\int_{t}^{T}\langle y(s)\rangle_{\epsilon }^{2k-1}\Big(\vert A(s) \vert\langle y(s) \rangle_{\epsilon }+\vert B(s) \vert \cdot \vert z(s) \vert  \notag \\
&&+\vert C(s) \vert \cdot \Vert \kappa(s) \Vert +\vert F(s,\alpha(s-))\vert\Big)\mathrm{d}s\bigg]  \notag \\
&\leq &\mathbb{E}[ \langle \xi \rangle_{\epsilon }^{2k}] +K_{0}\mathbb{E}\bigg [\int_{t}^{T}\Big(\langle y(s)\rangle_{\epsilon }^{2k}+\langle y(s) \rangle_{\epsilon }^{2k-1}\vert z(s)\vert  \notag \\
&&+\langle y(s) \rangle_{\epsilon }^{2k-1}\Vert\kappa(s) \Vert +\langle y(s)\rangle_{\epsilon }^{2k-1}\vert F(s,\alpha(s-))\vert\Big)\mathrm{d}s\bigg].\label{p1}
\end{eqnarray}
Here $K_{0}=K_{0}(K,k)$ is independent of $t$. Applying the
well-known Young's inequality ($ab\leq \frac{a^{p}}{p}+\frac{b^{q}}{q}$ with
$\frac{1}{p}+\frac{1}{q}=1$), we get
\begin{eqnarray*}
&&\mathbb{E}[ \langle y(t) \rangle_{\epsilon}^{2k}] +k(2k-1) \mathbb{E}\int_{t}^{T}\langle y(s) \rangle_{\epsilon }^{2k-2}\big[\vert z(s) \vert^{2}+\Vert \kappa(s)\Vert^{2}\big]\mathrm{d}s \\
&\leq &\mathbb{E}[ \langle \xi \rangle_{\epsilon }^{2k}] +K_{0}\mathbb{E}\bigg[\int_{t}^{T}\Big(\langle y(s)\rangle_{\epsilon }^{2k}+\frac{1}{\beta }\langle y(s) \rangle_{\epsilon }^{2k-2}\vert z(s)\vert^{2} \\
&&+\frac{1}{2}\langle y(s) \rangle_{\epsilon}^{2k-2}\Vert \kappa(s) \Vert^{2}+\vert F(s,\alpha(s-)) \vert^{2k}\Big)\mathrm{d}s\bigg],
\end{eqnarray*}
Let $\beta_{1}$ be large enough such that $K_{1}=k(2k-1) -\frac{K_{0}}{\beta_{1}}>0$. It immediately yields
\begin{eqnarray}
&&\mathbb{E}[ \langle y(t) \rangle_{\epsilon}^{2k}] +K_{1}\mathbb{E}\int_{t}^{T}\langle y(s)\rangle_{\epsilon }^{2k-2}\big[\vert z(s)\vert^{2}+\Vert \kappa(s) \Vert^{2}\big]\mathrm{d}s  \notag \\
&\leq &\mathbb{E}[ \langle \xi \rangle_{\epsilon }^{2k}] +K_{0}\mathbb{E}\bigg[\int_{t}^{T}\Big(\langle y(s)\rangle_{\epsilon }^{2k}+\vert F(s,\alpha(s-)) \vert^{2k}\Big)\mathrm{d}s\bigg].  \label{p2}
\end{eqnarray}
Similarly, let $\beta_{2}$ be large enough such that $\tilde{K}_{1}=k(2k-1) -\frac{K_{0}}{\beta_{2}}>0$. Then, we have
\begin{eqnarray}
&&\mathbb{E}\big[\langle y(t) \rangle_{\epsilon}^{2k}\big] +\tilde{K}_{1}\mathbb{E}\int_{t}^{T}\langle y(s) \rangle_{\epsilon }^{2k-2}\big[\vert z(s)\vert^{2}+\Vert \kappa(s)\Vert^{2}\big]\mathrm{d}s  \notag \\
&\leq &\mathbb{E}\big[\langle \xi \rangle_{\epsilon }^{2k}\big] +K_{0}\mathbb{E}\bigg[\int_{t}^{T}\Big(\langle y(s)\rangle_{\epsilon }^{2k}+\langle y(s) \rangle_{\epsilon }^{2k-1}\vert F(s,\alpha(s-))\vert\Big)\mathrm{d}s\bigg].  \label{p22}
\end{eqnarray}
Hence, it follows from Gronwall's inequality that
\begin{equation}
\mathbb{E}\big[\langle y(t) \rangle_{\epsilon}^{2k}\big] \leq K_{2}\bigg\{\mathbb{E}\Big(\langle \xi\rangle_{\epsilon }^{2k}\Big) +\mathbb{E}\bigg[\int_{0}^{T}\vert F(s,\alpha(s-))\vert^{2k}\mathrm{d}s\bigg]\bigg\},  \label{p3}
\end{equation}
and
\begin{eqnarray}
\mathbb{E}\bigg[\int_{t}^{T}\langle y(s) \rangle_{\epsilon }^{2k-2}\big[\vert z(s) \vert^{2}+\Vert \kappa(s) \Vert^{2}\big]\mathrm{d}s\bigg]
\leq K_{3}\bigg\{\mathbb{E}\big[\langle \xi \rangle_{\epsilon}^{2k}\big] +\mathbb{E}\bigg[\int_{0}^{T}\vert F(s,\alpha(s-)) \vert^{2k}\mathrm{d}s\bigg]\bigg\}.  \label{p33}
\end{eqnarray}
Note that $K_{2}=K_{2}(K,T,k) $. Besides since we postulate that
at the beginning being that $F$ are bounded, the above procedure becomes
valid (otherwise the integration on the right-hand side of (\ref{p2}) may
not exist; see (\ref{lp1}). Next, we want to refine the above estimate so
that(\ref{newest}) will follow. To this end, observe that (\ref{p3})
implies that its left-hand side is bounded uniformly in $t\in [ 0,T]$. Thus, it is allows us to define
\begin{equation}
\varphi(t) =\bigg\{ \sup_{0\leq s\leq t}\mathbb{E}\big[\langle y(s) \rangle_{\epsilon }^{2k}\big]\bigg\}^{\frac{1}{2k}},  \quad t\in [ 0,T].  \label{p4}
\end{equation}
We now come back (\ref{p22}), using (\ref{p4}). Define $\delta =\frac{1}{2K_{0}}-(1-\frac{1}{2k}) \tau $, for $\tau $ small enough. Then,
for any $t\in [ T-\delta,T]$, we have from Young's inequality, H\"{o}lder inequality and (\ref{p33})
\begin{eqnarray}
&&\varphi(t)^{2k}\leq \mathbb{E}\big[\langle \xi\rangle_{\epsilon }^{2k}\big] +K_{0}\mathbb{E}\bigg [t\varphi(t)^{2k} +\int_{t}^{T}\langle y(s) \rangle_{\epsilon}^{2k-1}\vert F(s,\alpha(s-)) \vert\mathrm{d}s\bigg ]  \notag \\
&\leq &\mathbb{E}\big[\langle \xi \rangle_{\epsilon }^{2k}\big] +\frac{1}{2}\varphi(t)^{2k}+K_{4}\bigg\{\int_{t}^{T}\Big(\mathbb{E}\vert F(s,\alpha(s-)) \vert^{2k}\Big)^{\frac{1}{2k}}\mathrm{d}s\bigg\}^{2k}.
\label{p5}
\end{eqnarray}
The constant $K_{4}=K_{4}(k,K,\delta) $ in (\ref{p5}) is independent of $t$. Hence, it follows from (\ref{p5}) that
\begin{equation}
\varphi(t)^{2k}\leq 2\mathbb{E}\big[\langle \xi\rangle_{\epsilon }^{2k}\big] +2K_{4}\bigg\{\int_{t}^{T}\Big(\mathbb{E}\vert F(s,\alpha(s-)) \vert^{2k}\Big)^{\frac{1}{2k}}\mathrm{d}s\bigg\}^{2k}.
\end{equation}
Now we repeat the same procedure on $[\delta,2\delta]$ and on $[2\delta,3\delta ] $, and so on. Eventually, we end up with
\begin{equation}
\varphi(t)^{2k}\leq \mathbb{E}\big[\langle \xi\rangle_{\epsilon }^{2k}\big] +K_{5}\bigg\{\int_{t}^{T}\Big(\mathbb{E}\vert F(s,\alpha(s-)) \vert^{2k}\Big)^{\frac{1}{2k}}\mathrm{d}s\bigg\}^{2k},
\end{equation}
with $K_{5}=K_{5}(k,\delta,K) $. According to the definition of (\ref{p4}), we conclude that
\begin{equation}
\sup_{0\leq s\leq t}\mathbb{E}\big[\langle y(s)\rangle_{\epsilon }^{2k}\big] \leq \mathbb{E}\big[\langle\xi \rangle_{\epsilon }^{2k}] +K_{5}\bigg\{ \int_{t}^{T}\Big(\mathbb{E}\vert F(s,\alpha(s-)) \vert^{2k}\Big)^{\frac{1}{2k}}\mathrm{d}s\bigg\}^{2k}.
\end{equation}
Letting $\epsilon \rightarrow 0$, we get the desired result. In the case
that we only have(\ref{lp1}), we may use the usual approach of
approximation. We thus complete the proof.
\hfill $\Box $

\noindent \textbf{Acknowledgements. }The author wish to thank Dr. Shuai Jing and Dr. Zhongyang Sun
for his conversation and suggestions.
%\paragraph{Proof of Lemma \protect\ref{estvs}}


\begin{thebibliography}{99}
\bibitem{Bis} J. M. Bismut, An introductory approach to duality in optimal
stochastic control. \emph{SIAM Rev.}, 20 (1978), 62-78.

\bibitem{CZE} Z. Chen and L Epstein, Ambiguity, risk, and asset returns in
continuous time. \emph{Econometrica}, 70 (2002), 1403-1443.

\bibitem{CIL} M. G. Crandall, H. Ishii, and P.-L. Lions, User's
guide to viscosity solutions of second order partial differential equations,
\emph{Bull. Amer. Math. Soc.}(N.S.),  27 (1992), 1-67.

\bibitem{CK} A. Cadenillas and I. Karatzas, The stochastic maximum principle
for linear, convex optimal control with random coefficients, \emph{SIAM J. Control Optim}., 33 (1995), 590-624.

\bibitem{CM} J. Cvitanic and J. Ma, Hedging options for a large investor and
forward-backward SDE's, \emph{The Annals of Applied Probability},  6 (1996), 370-398.

\bibitem{CC} D. Cuoco and J. Cvitanic, Optimal consumption choices for a
\textquotedblleft large\textquotedblright\ investor, \emph{Journal of
Economic Dynamics \& Control}, 22 (1998), 401-436.

\bibitem{DE} D. Duffie and L. G. Epstein, Stochastic differential utility.
\emph{Econometrica}. 60(1992) 353-394.

\bibitem{Don} C. Donnelly, Sufficient stochastic maximum principle in a
regime-switching diffusion model. \emph{Appl. Math. Optim.}, 64 (2011), 155-169.

\bibitem{DH} C. Donnelly and A. J. Heunis, Quadratic risk minimization in a
regime-switching model with portfolio constraints. \emph{SIAM J. Control Optim.}, 50 (2012), 2431-2461.

\bibitem{CE} S. N. Cohen, R.J. Elliott, Comparisons for backward stochastic
differential equations on Markov chains and related no-arbitrage conditions.
\emph{Ann. Appl. Probab.}, 20 (2010), 267-311.

\bibitem{EAM} R.J. Elliott, L. Aggoun, and J.B. Moore, Hidden Markov Models:
Estimation and Control, Springer, New York, 1994.

\bibitem{FS}
W.H. Fleming and H.M. Soner,
\textit{Controlled Markov Processes and Viscosity Solutions}, Springer, 2006.

\bibitem{HP} Y. Hu and S. Peng, Solution of forward-backward stochastic
differential equations. \emph{Probab. Theory Relat. Field}, 103 (1995), 273-283.

\bibitem{Hms} M. Hu, Stochastic global maximum principle for optimization
with recursive utilities, \emph{Probability, Uncertainty and Quantitative Risk}, (2017) 2:1.

\bibitem{HJX1} M. Hu, S. Ji, X. Xue, A Global Stochastic Maximum Principle
for Fully Coupled Forward-Backward Stochastic Systems, \emph{SIAM J. Control Optim.}, 56 (2018), 4309-4335.

\bibitem{HJX2} M. Hu, S. Ji, X. Xue, A note on the global stochastic maximum principle for
fully coupled forward-backward stochastic systems, arXiv:1812.10469.

\bibitem{JZ} S. Ji and X.Y. Zhou, A maximum principle for stochastic optimal
control with terminal state constrains, and its applications. \emph{Comm. Inf. Syst.}, 6 (2006), 321-337.

\bibitem{JX} S. Ji, X. Xue, The stochastic maximum principle in singular
optimal control with recursive utilities, \emph{J. Math. Anal. Appl.}, 471 (2019), 378-391.

\bibitem{KPQ} N.E. Karoui, S. Peng and M.C. Quenez, Backward stochastic
differential equations in nance. \emph{Math. Financ.}, 7 (1997) 1-71.

\bibitem{KZ} M. Kohlmann and X.Y. Zhou, Relationship between backward
stochastic differential equations and stochastic controls: A linear-quadratic approach.
\emph{SIAM J. Control Optim.}, 38 (2000) 1392-1407.

\bibitem{LZA} X. Li, X.Y. Zhou and M. Ait Rami,
Indefinite stochastic linear quadratic control with Markovian jumps in infinite time horizon.
\emph{Journal of Global Optimization}, 27 (2003), 149-175.

\bibitem{LZ} A. Lim and X.Y. Zhou, Linear-quadratic control of backward
stochastic differential equations. \emph{SIAM J. Control Optim.}, 40 (2001) 450-474.

\bibitem{LW} S. Lv and Z. Wu, Stochastic Maximum Principle for Forward-Backward Regime Switching Jump Diffusion Systems and
Applications to Finance. \emph{Chin. Ann. Math. Ser. B}, 39 (2018), 773-790.

\bibitem{Meng} Q. Meng, A maximum principle for optimal control problem of
fully coupled forward-backward stochastic systems with partial information,
\emph{Science in China. Series A}, 52 (2009), 1579-1588.

\bibitem{PP1} E. Pardoux and S. Peng, Adapted solution of a backward
stochastic differential equation. \emph{Syst. Control Lett.}, 14 (1990), 55-61.

\bibitem{Pa} O. M. Pamen,  Maximum principles of Markov regime-switching forward-backward stochastic differential equations with jumps and partial information. \emph{J. Optim Theory Appl.},  175 (2017), 373-410.

\bibitem{PZ} S. Peng and Z. Wu, Fully coupled forward-backward stochastic dierential equations and applications to optimal control. \emph{SIAM J.
Control Optim.}, 37(1999) 825-843.

\bibitem{P1} S. Peng, Backward stochastic differential equations and
applications to optimal control. \emph{Appl. Math. Optim.}, 27(1993), 125-144.

\bibitem{QZ} Q. Zhang, Stock trading: an optimal selling rule. \emph{SIAM J.
Control Optim.}, 40 (2001) 64-87.

\bibitem{RH} R. Buckdahn and Y. Hu, Hedging contingent claims for a large investor in an incomplete market,
\emph{Advances in Applied Probability}, 30 (1998), 239-255.

\bibitem{SW} J. Shi and Z. Wu, The maximum principle for fully coupled
forward-backward stochastic control system. \emph{Acta Automat. Sinica}, 32 (2006), 161-169.

\bibitem{SKM} Z. Sun, I. Kemajou-Brown and O. Menoukeu-Pamen,
A risk-sensitive maximum principle for a Markov regime-switching jump-di usion system and applications. \emph{ESAIM-Control Optim. Calc. Var.},24 (2018), 985-1013.

\bibitem{SGZ} Z. Sun, J. Guo and X. Zhang, Maximum Principle for Markov
Regime-Switching Forward--Backward Stochastic Control System with Jumps and
Relation to Dynamic Programming, \emph{J Optim Theory Appl.}, 176 (2018), 319-350.

\bibitem{TW} R. Tao and Z. Wu, Maximum principle for optimal control
problems of forward-backward regime-switching system and applications. \emph{Syst. Control Lett.}, 61 (2012), 911-917.

\bibitem{ZXCS} D. Zhu, Y. Xie, W. K. Ching and T.K Siu, Optimal portfolios
with maximum Value-at-Risk constraint under a hidden Markovian
regime-switching model. \emph{Automatica}, 74 (2016), 194-205.

\bibitem{Zhu} C. Zhu, Optimal control of the risk process in a
regime-switching environment. \emph{Automatica}, 47 (2011), 1570-1579.

\bibitem{Peng1990} S. Peng, A general stochastic maximum principle for
optimal control problems. \emph{SIAM J. Control Optim.}, 28 (1990), 966-979.

\bibitem{PW} S. Peng, Z. Wu, Fully Coupled Forward-Backward Stochastic
Differential Equations and Applications to Optimal Control. \emph{SIAM J. Control
Optim.}, 37 (1999), 825-843.

\bibitem{SZ} Q. Song, C. Zhu, On singular control problems with state
constraints and regime-switching. \emph{Automatica}, 70 (2016), 66-73.

\bibitem{Wu2013} Z. Wu, A general maximum principle for optimal control of
forward backward stochastic systems. \emph{Automatica}, 49 (2013),1473-1480.

\bibitem{Wu1998} Z. Wu, Maximum principle for optimal control problem of
fully coupled forward-backward stochastic systems. \emph{Syst. Sci. Math.Sci.}, 11 (1998), 249-259.

\bibitem{Xu} W. Xu, Stochastic maximum principle for optimal control problem
of forward and backward system. \emph{J. Austral. Math. Soc. Ser. B} 37 (1995), 172-185.

\bibitem{YZ} J. Yong and X. Zhou, Stochastic controls: Hamiltonian systems
and HJB equations. Vol. 43. Springer-Verlag, New York, 1999.

\bibitem{Yong2010} J. Yong, Optimality variational principle for controlled
forward-backward stochastic differential equations with mixed
initial-terminal conditions, \emph{SIAM J. Control Optim.}, 48 (2010), 4119-4156.

\bibitem{ZEK} X. Zhang, R. J. Elliott, T. K. Siua, Stochastic maximum
principle for a Markov regime-switching jump-diffusion model and its
application to finance, \emph{SIAM J. Control Optim.}, 50 (2012), 964-990.

%\bibitem{ZXCS} D. Zhu, Y. Xie, W. K. Ching and T.K Siu, Optimal portfolios
%with maximum Value-at-Risk constraint under a hidden Markovian
%regime-switching model. \emph{Automatica}. 74(2016) 194-205.

\bibitem{ZC} C. Zhu, Optimal control of the risk process in a
regime-switching environment. \emph{Automatica}, 47 (2011) ,1570-1579.

\bibitem{ZSM} X. Zhang, T.K. Siu and Q. B. Meng. Portfolio selection in the
enlarged Markovian regime-switching market. \emph{SIAM J. Control Optim.}, 48 (2010), 3368-3388.

\bibitem{ZEG} X. Zhang, R.J. Elliott, T.K. Siu and J.Y. Guo. Markovian
regime-switching market completion using additional markov jump assets.
\emph{IMA Journal of Management Mathematics}, 23 (2011), 283-305.

\bibitem{ZSX} X. Zhang, Z. Sun, and J. Xiong. A general stochastic maximum
principle for a Markov regime switching jump-diffusion model of mean-field
type. \emph{SIAM J. Control Optim.}, 56 (2018), 2563-2592.
\end{thebibliography}
\end{document}